# EFFICIENT PREDICTION FOR LINEAR AND NONLINEAR AUTOREGRESSIVE MODELS

By Ursula U. Müller, Anton Schick[1] and Wolfgang Wefelmeyer

*Texas A&M University, Binghamton University and Universität zu Köln*

Conditional expectations given past observations in stationary time series are usually estimated directly by kernel estimators, or by plugging in kernel estimators for transition densities. We show that, for linear and nonlinear autoregressive models driven by independent innovations, appropriate smoothed and weighted von Mises statistics of residuals estimate conditional expectations at better parametric rates and are asymptotically efficient. The proof is based on a uniform stochastic expansion for smoothed and weighted von Mises processes of residuals. We consider, in particular, estimation of conditional distribution functions and of conditional quantile functions.

**1. Introduction.** Let $X_0, \dots, X_n$ be observations from a real-valued stationary time series. Conditional expectations $E(q(X_{n+m})|X_n = x)$ with lag $m$ of some known function $q$ can be estimated by kernel estimators. For asymptotic results under various mixing conditions, we refer to [8, 23, 24, 29, 36, 37, 43, 44]. If the time series is first-order Markov with transition density $p(x,y)$, a conditional expectation of $q$ with lag one can be written $E(q(X_{n+1})|X_n = x) = \int q(y)p(x,y)\,dy$ and it can be estimated by plugging in a kernel estimator $\hat{p}(x,y)$. Asymptotic results for such estimators of conditional expectations are in [9, 16, 26, 27, 28].

If the Markov chain follows a nonparametric autoregressive model $X_i = r(X_{i-1}) + \varepsilon_i$, with unknown autoregression function $r$ and independent and identically distributed (i.i.d.) mean zero innovations $\varepsilon_i$, then $E(q(X_{n+1})|X_n =$

Received May 2003; revised April 2005.
[1]Supported in part by NSF Grants DMS-00-72174 and DMS-04-05791.

$x) = E[q(\varepsilon_1 + r(x))]$. Let $\tilde r$ denote a (kernel) estimator of the autoregression function. Write $\tilde\varepsilon_i = X_i - \tilde r(X_{i-1})$ for the residuals and $\tilde F(y) = \frac{1}{n}\sum_{i=1}^n \mathbf{1}[\tilde\varepsilon_i \le y]$ for the empirical distribution function based on them. The representation suggests estimating the conditional expectation by an empirical estimator

$$\frac{1}{n}\sum_{i=1}^n q(\tilde\varepsilon_i + \tilde r(x)) = \int q(y + \tilde r(x))\, d\tilde F(y). \tag{1.1}$$

The convergence rate of (1.1) is given by the convergence rate of $\tilde r$.

Suppose now that we have a linear or nonlinear *parametric* model $r = r_\vartheta$ for the autoregression function. In this case we can use a $n^{1/2}$-consistent estimator $\tilde\vartheta$ for $\vartheta$ and the $n^{1/2}$-consistent estimator $\tilde r = r_{\tilde\vartheta}$ for $r$, and we can estimate the innovations $\varepsilon_i$ by $\tilde\varepsilon_i = X_i - r_{\tilde\vartheta}(X_{i-1})$. Under appropriate smoothness and integrability conditions on the function $q$, one can prove by a Taylor expansion that the resulting estimator (1.1) is $n^{1/2}$-consistent; see [33] for closely related details in a different problem. In particular, the estimator (1.1) converges at a faster rate than the nonparametric estimators. If $\tilde\vartheta$ is asymptotically normal, so is (1.1). Such results could also be obtained for heteroscedastic autoregressive models $X_i = r_\vartheta(X_{i-1}) + s_\vartheta(X_{i-1})\varepsilon_i$, including ARCH models, and for GARCH models. For GARCH models and smooth $q$, one could use limit results for the empirical process of residuals obtained by Boldin [6, 7] and Berkes and Horváth [2, 3].

Since the innovations are assumed to have mean zero, the residual-based empirical distribution function $\tilde F$ is not an efficient estimator of $F$. Thus, improvements over (1.1) are possible by replacing $\tilde F$ by an efficient estimator. Here efficiency is meant in the sense of a semiparametric version of Hájek and Le Cam's convolution theorem; see also Section 6. An efficient estimator of $F$ has been constructed in [32], but this estimator is not a distribution function. Alternative efficient estimators that are distribution functions are discussed in [18]. One such estimator is the weighted residual-based empirical distribution function

$$\tilde F_w(y) = \frac{1}{n}\sum_{i=1}^n w_i \mathbf{1}[\tilde\varepsilon_i \le y], \qquad y \in \mathbb{R},$$

with an efficient estimator $\tilde\vartheta$ and random weights $w_i$ chosen following the empirical likelihood approach of Owen [19, 20] so that, with probability tending to one, $\tilde F_w$ has mean zero, that is, $\int y\, d\tilde F_w(y) = (1/n)\sum_{i=1}^n w_i\tilde\varepsilon_i = 0$. The resulting weighted version of (1.1) is the estimator

$$\int q(y + r_{\tilde\vartheta}(x))\, d\tilde F_w(y) = \frac{1}{n}\sum_{i=1}^n w_i q(\tilde\varepsilon_i + r_{\tilde\vartheta}(x)). \tag{1.2}$$

This estimator is efficient if $\tilde\vartheta$ is. This is a consequence of the fact that smooth functionals of efficient estimators are efficient. An alternative to weighting would be to subtract an appropriate "estimator of zero" from the estimator that corrects the influence function. See [15] and [12] for models with i.i.d. data; [17] for Markov chains; and [32, 33] for time series residuals. However, weighting has the advantage that, with high probability, the information of mean zero is used exactly, so we expect better small-sample properties.

Let us now look at some special cases in which simple alternative estimators are also available. For the conditional mean of lag one, for which $q(x) = x$, we have $E(X_{n+1}|X_n = x) = r_\vartheta(x)$. This can be estimated directly by $r_{\tilde\vartheta}(x)$. This estimator is efficient if $\tilde\vartheta$ is. The estimator (1.1) is $(1/n)\sum_{i=1}^n \tilde\varepsilon_i + r_{\tilde\vartheta}(x)$, which is $n^{1/2}$-consistent but is not efficient even if $\tilde\vartheta$ is. The weighted estimator (1.2) equals the direct estimator $r_{\tilde\vartheta}(x)$ with probability tending to one. Hence, it is efficient if $\tilde\vartheta$ is. Another special case is the conditional second moment of lag one, for which $q(x) = x^2$. We have $E(X_{n+1}^2|X_n = x) = E[\varepsilon_1^2] + r_\vartheta^2(x)$. The empirical estimator (1.1) is $(1/n)\sum_{i=1}^n(\tilde\varepsilon_i + r_{\tilde\vartheta}(x))^2$. It is $n^{1/2}$-consistent, but not efficient. A more direct $n^{1/2}$-consistent estimator is the plug-in estimator $(1/n)\sum_{i=1}^n \tilde\varepsilon_i^2 + r_{\tilde\vartheta}^2(x)$. However, it is not efficient in general even if $\tilde\vartheta$ is, since it does not (fully) exploit the fact that the innovations have mean zero. Efficient estimators are given by the weighted empirical estimator and the (asymptotically equivalent) weighted plug-in estimator $(1/n)\sum_{i=1}^n w_i\tilde\varepsilon_i^2 + r_{\tilde\vartheta}^2(x)$, both with efficient $\tilde\vartheta$.

Similar results are possible for lag two. The conditional expectation $E(q(X_{n+2})|X_n = x)$ becomes $E[q(\varepsilon_2 + r_\vartheta(\varepsilon_1 + r_\vartheta(x)))]$ and can be estimated $n^{1/2}$-consistently by the von Mises statistic

$$\iint q(z + r_{\tilde\vartheta}(y + r_{\tilde\vartheta}(x)))\, d\tilde F(y) d\tilde F(z) = \frac{1}{n^2}\sum_{i=1}^n\sum_{j=1}^n q(\tilde\varepsilon_j + r_{\tilde\vartheta}(\tilde\varepsilon_i + r_{\tilde\vartheta}(x)))$$

and the weighted von Mises statistic

$$\iint q(z + r_{\tilde\vartheta}(y + r_{\tilde\vartheta}(x)))\, d\tilde F_w(y)\, d\tilde F_w(z)$$
$$= \frac{1}{n^2}\sum_{i=1}^n\sum_{j=1}^n w_i w_j q(\tilde\varepsilon_j + r_{\tilde\vartheta}(\tilde\varepsilon_i + r_{\tilde\vartheta}(x))).$$

The latter will be efficient if an efficient estimator $\tilde\vartheta$ of $\vartheta$ is used. The von Mises statistics are easier to use than the usual kernel estimator because they do not require a choice of bandwidth. For certain $q$, simpler alternative estimators are available. For example, the conditional mean of lag two equals $E[r_\vartheta(\varepsilon_1 + r_\vartheta(\mathbf{x}))]$ and can be estimated more directly by the average $(1/n)\sum_{i=1}^n r_{\tilde\vartheta}(\tilde\varepsilon_i + r_{\tilde\vartheta}(\mathbf{x}))$ or the weighted average $(1/n)\sum_{i=1}^n w_i r_{\tilde\vartheta}(\tilde\varepsilon_i +$ $r_{\tilde\vartheta}(\mathbf{x})$). The latter coincides with the weighted von Mises statistic with probability tending to one. A degenerate case would be the linear AR(1) model, with $r_\vartheta(x) = \vartheta x$, for which the conditional mean of lag two is $\vartheta^2 x$, which is estimated efficiently by $\tilde\vartheta^2 x$ with $\tilde\vartheta$ efficient for $\vartheta$. The weighted von Mises statistic coincides with this simple efficient estimator with probability tending to one. Simplified versions of the von Mises statistics are also available for estimating higher conditional moments of lag two. The conditional second moment of lag two simplifies to $E[\varepsilon_1^2] + E[r_\vartheta^2(\varepsilon_1 + r_\vartheta(x))]$ and can be estimated $n^{1/2}$-consistently by the average $(1/n)\sum_{i=1}^n(\tilde\varepsilon_i^2 + r_{\tilde\vartheta}^2(\tilde\varepsilon_i + r_{\tilde\vartheta}(x)))$ or the weighted average $(1/n)\sum_{i=1}^n w_i(\tilde\varepsilon_i^2 + r_{\tilde\vartheta}^2(\tilde\varepsilon_i + r_{\tilde\vartheta}(x)))$. The latter equals the weighted von Mises estimator with probability tending to one and is efficient if $\tilde\vartheta$ is.

The above shows that conditional expectations of lags one and two can be estimated $n^{1/2}$-consistently and efficiently for smooth $q$ in nonlinear autoregression models of order one. To prove $n^{1/2}$-consistency of the estimator (1.1) for more general $q$, we need an appropriate balance of smoothness assumptions on $q$ and on the innovation distribution. For discontinuous $q$, we must assume that the innovations have a smooth density $f$. One may then also want to replace $\tilde F$ and $\tilde F_w$ by smoothed versions $\tilde F_s$ and $\tilde F_{sw}$, say, $d\tilde F_s(y) = \tilde f(y)\,dy$ and $d\tilde F_{sw}(y) = \tilde f_w(y)\,dy$, where $\tilde f$ is a kernel estimator $\tilde f(y) = (1/n)\sum_{i=1}^n k_{b_n}(y - \tilde\varepsilon_i)$ of the density $f$ and $\tilde f_w$ is a weighted kernel estimator $\tilde f_w(y) = (1/n)\sum_{i=1}^n w_i k_{b_n}(y - \tilde\varepsilon_i)$. Here $k_{b_n}(y) = k(y/b_n)/b_n$ for some kernel $k$ and some bandwidth $b_n$. These kernel estimators were studied in [18]. Efficiency of the smoothed and weighted residual-based empirical distribution function $\tilde F_{sw}$ was also shown there. The resulting smoothed and weighted von Mises statistic

$$\iint q(z + r_{\tilde\vartheta}(y + r_{\tilde\vartheta}(x)))\,\tilde f_w(y)\,dy\,\tilde f_w(z)\,dz$$

preserves $n^{1/2}$-consistency and efficiency even though the kernel estimators have a slower rate of convergence. Simulations show that smoothing improves the small-sample behavior of our estimator noticeably, especially if $q$ is not smooth (see Table 1). This is a second-order effect. For theoretical results in this direction, see [11]. We note that the choice of bandwidth is less critical here than for the usual kernel estimators. In particular, the asymptotic variance of our estimator does not depend on the choice of bandwidth in the allowed range.

The smoothed and weighted estimator

$$\iint q(z + r_{\tilde\vartheta}(y + r_{\tilde\vartheta}(x)))\,\tilde f_w(y)\,dy\,\tilde f_w(z)\,dz$$

equals

$$\frac{1}{n^2}\sum_{i=1}^n\sum_{j=1}^n w_i w_j \iint q(\tilde\varepsilon_j + b_n u + r_{\tilde\vartheta}(\tilde\varepsilon_i + b_n v + r_{\tilde\vartheta}(x)))k(u)\,du\,k(v)\,dv.$$

TABLE 1
*Simulated mean squared error for various von Mises estimators*

| | n | U | W | 1.50 | 1.75 | 2.00 | 2.25 | 2.50 | 2.75 |
|---|---|---|---|---|---|---|---|---|---|
| Normal | 50 | 6181 | 967 | 512 | 462 | 430 | 414 | 411 | 417 |
| | 100 | 3153 | 460 | 299 | 279 | 266 | 261 | 264 | 273 |
| | 200 | 1615 | 227 | 168 | 160 | 156 | 155 | 160 | 168 |
| Logistic | 50 | 6184 | 1218 | 647 | 591 | 558 | 544 | 545 | 558 |
| | 100 | 3204 | 606 | 390 | 367 | 356 | 355 | 364 | 380 |
| | 200 | 1620 | 296 | 220 | 213 | 212 | 217 | 227 | 243 |
| $t(5)$ | 50 | 6363 | 1513 | 803 | 738 | 701 | 686 | 690 | 706 |
| | 100 | 3234 | 756 | 495 | 470 | 459 | 461 | 474 | 495 |
| | 200 | 1646 | 375 | 281 | 274 | 275 | 283 | 299 | 320 |

The table entries are $10^6 \times$MSE of the von Mises estimator (U), the weighted von Mises estimator (W) and the smoothed and weighted von Mises estimator for different bandwidths $b_n = cn^{-1/4}$ with $c = 1.5, 1.75, 2, 2.25, 2.5, 2.75$. The simulations are based on 20,000 repetitions. We estimate the conditional probability $P(X_{n+2} \le 0 | X_n = 0.5)$ in the AR(1) model $X_i = \vartheta X_{i-1} + \varepsilon_i$ with $\vartheta = 0.5$ for sample sizes $n = 50, 100, 200$. The innovation distributions are the standard normal distribution, the logistic distribution and the $t$-distribution with five degrees of freedom, the latter two scaled to have variance one. As estimator of $\vartheta$, the sample autocorrelation coefficient, was used. The standard error of a simulated MSE is about 1% of the MSE.

When the latter double integral is difficult to calculate, it can be approximated by Riemann sums, resulting in

$$\frac{4}{(nN)^2} \sum_{i=1}^{n} \sum_{j=1}^{n} \sum_{s=1}^{N} \sum_{t=1}^{N} w_i w_j q(\tilde{\varepsilon}_j + b_n u_s + r_{\tilde{\vartheta}}(\tilde{\varepsilon}_i + b_n u_t + r_{\tilde{\vartheta}}(x))) k(u_s) k(u_t).$$

Here $u_1, \dots, u_N$ denote the midpoints of a partition of the compact support $[-1, 1]$ of the kernel $k$ into $N$ intervals of equal lengths. This shows that the smoothed estimator is easy to compute.

Weighting can lead to drastic variance reductions, especially if $q$ is asymmetric, for example, for odd moments and for distribution functions. See Example 3.2 and Example 5.5, which treat smoothed and weighted von Mises statistics in the classical autoregressive model of order one. Example 3.2 reports a possible variance reduction of up to 64% for the one-lag conditional distribution function. Similar improvements through weighting are obtained for estimators of expectations under the innovation distribution; see [18], Sections 4 and 5. Example 5.5 shows that variance reductions of over 98% are possible in the case of estimating the lag-two conditional distribution function. The simulation results in Table 1 show that, for small to moderate sample sizes, the actual variance reductions might be even larger due to the second-order effect of smoothing.

It is the purpose of this paper to extend and sharpen the results on smoothed and weighted von Mises statistics outlined above in several directions: to linear and nonlinear autoregressive models of higher order, to conditional expectations with higher lags, to functions $q$ of more than one argument and to uniform results over classes of functions. We are particularly interested in estimating univariate and multivariate conditional distribution functions. They give rise to $n^{1/2}$-consistent estimators of conditional quantiles. Other applications are conditional probabilities of staying in a certain band, for example, $P(|X_{n+1}-x|\le c_1, |X_{n+2}-x|\le c_2|X_n=x)$, or conditional probabilities that the time series increases over a certain period, for example, $P(X_{n+3}>X_{n+2}>X_{n+1}>x|X_n=x)$.

Specifically, we consider linear or nonlinear autoregressive models of order $p$,

$$X_i = r_\vartheta(\mathbf{X}_{i-1}) + \varepsilon_i, \tag{1.3}$$

with $\mathbf{X}_{i-1}=(X_{i-p},\dots,X_{i-1})$ and $\vartheta$ a $d$-dimensional parameter, and we construct estimators for conditional expectations $E(q(X_{n+1},\dots,X_{n+m})|\mathbf{X}_n = \mathbf{x})$ for some known function $q$ of $m$ arguments and some fixed vector $\mathbf{x} = (x_1,\dots,x_p)$. Using the representation of the autoregressive process, such conditional expectations can be written

$$E(q(X_{n+1},\dots,X_{n+m})|\mathbf{X}_n=\mathbf{x}) = E[q(\varrho_\vartheta(\varepsilon_{n+1},\dots,\varepsilon_{n+m}))]$$

for some function $\varrho_\vartheta$. For lag two, that is, $m=2$, we have $\varrho_\vartheta(\varepsilon_1,\varepsilon_2) = (\varepsilon_1 + r_\vartheta(\mathbf{x}), \varepsilon_2 + r_\vartheta(x_2,\dots,x_p,\varepsilon_1+r_\vartheta(\mathbf{x})))$. Let $\tilde\vartheta$ be a $n^{1/2}$-consistent estimator of $\vartheta$. Using it, we can form the residuals $\tilde\varepsilon_i = X_i - r_{\tilde\vartheta}(\mathbf{X}_{i-1})$, $i=1,\dots,n$. We estimate the conditional expectations by the smoothed and weighted von Mises statistic

$$\int\cdots\int q(\varrho_{\tilde\vartheta}(y_1,\dots,y_m))\prod_{j=1}^m \tilde f_w(y_j)\,dy_j.$$

It is efficient if an efficient estimator $\tilde\vartheta$ for $\vartheta$ is used. We obtain $n^{1/2}$-consistency and asymptotic normality not just for fixed $q$, but uniformly over large classes of functions. We show, in particular, that our estimator, viewed as a stochastic process indexed by $q$ and suitably standardized, converges to a Gaussian process. This is in contrast to the usual kernel estimators, for which limit theorems can hold only locally, in intervals shrinking in proportion to the bandwidth $b_n$.

Independence of innovations has recently also been exploited for other functionals. Schick and Wefelmeyer [33] use this idea to reduce the variance in estimating linear functionals of the stationary law of invertible linear processes. Saavedra and Cao [31] obtain a $n^{1/2}$-consistent estimator for the stationary density of an MA(1) process. Schick and Wefelmeyer [34] prove asymptotic efficiency of a modified version and Schick and Wefelmeyer [35] obtain functional central limit theorems in the case of MA($q$) processes, considering the density as an element of the function space $L_1$ or $C_0$. For nonparametric regression, van Keilegom and Veraverbeke [41, 42] and Van Keilegom, Akritas and Veraverbeke [40] exploit independence of the error and the covariate to obtain improved estimators for the conditional density, distribution function and hazard rate of the response given the covariate.

The paper is organized as In Section 2 we derive a stochastic expansion for smoothed von Mises processes based on residuals,

$$\psi(h,\tilde f)=\int\cdots\int h(y_1,\dots,y_m)\prod_{j=1}^m \tilde f(y_j)\,dy_j,$$

and for weighted versions $\psi(h,\tilde f_w)$, uniform over appropriate classes $\mathcal{H}$ of functions $h$. These are results of independent interest. To describe them, let $\hat f(y)=\frac{1}{n}\sum_{i=1}^n k_{b_n}(y-\varepsilon_i)$ be the kernel estimator based on the actual innovations, and

$$\bar h(y)=E(h(\varepsilon_1,\dots,\varepsilon_m)|\varepsilon_1=y)+\cdots+E(h(\varepsilon_1,\dots,\varepsilon_m)|\varepsilon_m=y).$$

The expansion of $\psi(h,\tilde f)$ is of the form

$$\psi(h,\tilde f)-\psi(h,f)=\int \bar h(y)(\hat f(y)-f(y))\,dy+D(h)^\top(\tilde\vartheta-\vartheta)+R_n(h),$$

with $\sup_{h\in\mathcal{H}}|R_n(h)|=o_p(n^{-1/2})$. Here $D(h)=E[\bar h(\varepsilon)\ell(\varepsilon)]E[\dot r_\vartheta(\mathbf{X})]$, where $\ell=-f'/f$ is the score function for location of the innovation distribution, $\dot r_\vartheta(\mathbf{X})$ is the gradient of $r_\vartheta(\mathbf{X})$ with respect to $\vartheta$ and $(\mathbf{X},\varepsilon)$ is short for $(\mathbf{X}_0,\varepsilon_1)$. The expansion of the weighted version differs as

$$\psi(h,\tilde f_w)=\psi(h,\tilde f)-\frac{E[\varepsilon\bar h(\varepsilon)]}{\sigma^2}\left(\frac{1}{n}\sum_{i=1}^n\varepsilon_i-E[\dot r_\vartheta(\mathbf{X})]^\top(\tilde\vartheta-\vartheta)\right)+R_{nw}(h),$$

where again $\sup_{h\in\mathcal{H}}|R_{nw}(h)|=o_p(n^{-1/2})$. In the above expansions, the terms involving $\tilde\vartheta-\vartheta$ come from replacing the estimated innovations by the true ones. Note that $\{n^{1/2}\int\bar h(y)(\hat f(y)-f(y))\,dy : h\in\mathcal{H}\}$ is a smoothed empirical process. Such processes have been studied by Yukich [46], van der Vaart [38], Rost [25] and Radulović and Wegkamp [21, 22]. They give conditions under which the smoothed empirical process is asymptotically equivalent to the usual empirical process. We refer to the book by van der Vaart and Wellner [39] for a general overview of empirical processes. We have an envelope and Lebesgue densities and give, in Propositions 2.1 and 2.2, versions of the results of van der Vaart [38] and Rost [25] with simpler assumptions. Together with the above expansions, these results imply that if $\tilde\vartheta$ is asymptotically linear, then so are the von Mises process $\{n^{1/2}(\psi(h,\tilde f)-\psi(h,f)) : h\in$ $\mathcal{H}\}$ and its weighted version $\{n^{1/2}(\psi(h,\tilde{f}_w)-\psi(h,f)) : h\in\mathcal{H}\}$. This implies that these processes converge weakly to tight Gaussian processes. For our applications to estimation of conditional expectations, we need versions in which the function $h$ is indexed by $\vartheta$ and $q$. We formulate such results in Theorem 2.2.

In Sections 3 to 5 we apply our results on von Mises processes to estimation of conditional expectations of lags one and two. We get by with mild assumptions on the innovation density and the autoregression function. In particular, we cover discontinuous autoregression functions such as those appearing in self-exciting threshold autoregressive (SETAR) models. Higher lags can be treated along these lines, but the stochastic expansions of the estimators are notationally cumbersome. In particular, Theorem 3.1 specializes Theorem 2.2 to the case of estimating conditional expectations of lag one. In Theorems 3.2 and 3.3 we apply Theorem 3.1 to estimators for conditional distribution functions and for the conditional expectation of a fixed function $q$. Theorems 4.1, 5.1 and 5.2 give analogous results for conditional expectations and conditional distribution functions of lag two. Examples 3.1 and 5.4 apply these results to conditional quantile processes of lags one and two. Our results are new and nontrivial, even for the linear autoregressive model of order one.

In Section 6 we show that the weighted versions of our estimators are efficient if an efficient estimator for $\vartheta$ is used. This is done by checking that the influence function then equals the efficient influence function for estimating $\psi(h_\vartheta, f)$ with $h_\vartheta = q\circ\varrho_\vartheta$. Efficient estimators for $\vartheta$ in nonlinear autoregression with mean zero innovations are constructed in [14].

Section 7 contains two technical lemmas. Lemma 7.1 gives a characterization of compact subsets of $L_2(\nu)$ for measures $\nu$ with Lebesgue density. It says that a closed subset of $L_2(\nu)$ with an envelope translation-continuous at zero is compact if and only if the subset is equi-translation-continuous at zero. Lemma 7.2 gives conditions for uniform differentiability of integrals with respect to Hellinger differentiable densities.

**2. Smoothed and weighted von Mises processes of residuals.** Consider observations $X_{1-p},\dots,X_n$ from a stationary and ergodic nonlinear autoregressive process $X_i = r_\vartheta(\mathbf{X}_{i-1})+\varepsilon_i$ of order $p$, where $\mathbf{X}_{i-1}=(X_{i-p},\dots,X_{i-1})$ and $\vartheta$ is a $d$-dimensional parameter. Assume that the innovations $\varepsilon_i$ are i.i.d. with mean zero, finite variance $\sigma^2$ and positive density $f$ and are independent of $\mathbf{X}_0$. Let $\tilde{\vartheta}$ be a $n^{1/2}$-consistent estimator for $\vartheta$. Estimate the innovations $\varepsilon_i$ by residuals $\tilde{\varepsilon}_i = X_i - r_{\tilde{\vartheta}}(\mathbf{X}_{i-1})$ and the innovation density $f$ by the kernel estimator

$$\tilde{f}(y)=\frac{1}{n}\sum_{i=1}^{n}k_{b_n}(y-\tilde{\varepsilon}_i)$$

or the weighted kernel estimator

$$\tilde{f}_w(y) = \frac{1}{n}\sum_{i=1}^{n} w_i k_{b_n}(y - \tilde{\varepsilon}_i),$$

where $k_{b_n}(y) = k(y/b_n)/b_n$ for a kernel $k$ and a bandwidth $b_n$. Following Owen [19, 20], we choose positive weights $w_i$ of the form

$$w_i = \frac{1}{1 + \tilde{\lambda}\tilde{\varepsilon}_i},$$

where $\tilde{\lambda}$ is chosen such that $\sum_{i=1}^{n} w_i \tilde{\varepsilon}_i = 0$. By Müller, Schick and Wefelmeyer [18], this is possible with probability tending to one. When there is no solution, we set $\tilde{\lambda} = 0$.

In this section we obtain a uniform stochastic expansion for smoothed von Mises processes based on residuals $\tilde{\varepsilon}_1, \dots, \tilde{\varepsilon}_n$,

$$\psi(h, \tilde{f}) = \int \cdots \int h(y_1, \dots, y_m) \prod_{j=1}^{m} \tilde{f}(y_j)\, dy_j,$$

and their weighted versions $\psi(h, \tilde{f}_w)$. Here the index $h$ runs through a family $\mathcal{H}$ of functions from $\mathbb{R}^m$ to $\mathbb{R}$ with envelope $H$, that is, $|h| \le H$ for all $h \in \mathcal{H}$. We assume the envelope to be of the form

$$H(y_1, \dots, y_m) = V(y_1) \cdots V(y_m), \tag{2.1}$$

where $V$ is a measurable function satisfying the following conditions.

Assumption V. The function $V$ satisfies $V \ge 1$ and, for some $\alpha > 1$,

$$\int (1 + |y|)^{\alpha} V^2(y) f(y)\, dy < \infty.$$

Moreover, the function $D$ defined by

$$D(s) := \sup_{y \in \mathbb{R}} \frac{|V(y+s) - V(y)|}{V(y)}, \qquad s \in \mathbb{R},$$

is bounded on compacts and is continuous at 0,

$$D(s) \to 0 \qquad \text{as } s \to 0. \tag{2.2}$$

If $V = 1$, then Assumption V is satisfied with $\alpha = 2$. Another example of a function satisfying Assumption V is $V(y) = (1 + |y|)^{\gamma}$ with $\gamma \ge 0$, provided $\int |y|^{2\gamma + \alpha} f(y)\, dy$ is finite for some $\alpha > 1$.

Write $\varepsilon$ and $\mathbf{X}$ for random variables with the same joint distribution as $\varepsilon_i$ and $\mathbf{X}_{i-1}$. Denote the distribution functions of $\varepsilon$ and $\mathbf{X}$ by $F$ and $G$. We make the following assumptions on the density $f$ and the autoregression function $r_{\vartheta}$.

ASSUMPTION F. The density $f$ has finite Fisher information for location, that is, $f$ is absolutely continuous with almost everywhere derivative $f'$, and $E[\ell^2(\varepsilon)] = \int \ell^2 \, dF$ is finite, where $\ell = -f'/f$.

ASSUMPTION R. The function $\tau \mapsto r_\tau(\mathbf{x})$ is continuously differentiable for all $\mathbf{x}$ with gradient $\tau \mapsto \dot r_\tau(\mathbf{x})$. For each constant $C$,

$$\begin{aligned} &\sup_{|\tau-\vartheta|\le Cn^{-1/2}} \sum_{i=1}^n (r_\tau(\mathbf{X}_{i-1}) - r_\vartheta(\mathbf{X}_{i-1}) - \dot r_\vartheta(\mathbf{X}_{i-1})^\top(\tau-\vartheta))^2 \\ &\qquad = O_p(n^{-2/3}). \end{aligned} \tag{2.3}$$

Moreover, $E[|\dot r_\vartheta(\mathbf{X})|^{5/2}] = \int |\dot r_\vartheta|^{5/2} \, dG < \infty$ and the matrix $E[\dot r_\vartheta(\mathbf{X}) \dot r_\vartheta(\mathbf{X})^\top] = \int \dot r_\vartheta \dot r_\vartheta^\top \, dG$ is positive definite.

A sufficient condition for (2.3) is a Hölder condition with exponent 2/3 on the gradient $\dot r_\tau$,

$$|\dot r_\tau(\mathbf{x}) - \dot r_\vartheta(\mathbf{x})| \le |\tau - \vartheta|^{2/3} A(\mathbf{x}),$$

with $A \in L_2(G)$.

Finally, we impose the following assumptions on the kernel and the bandwidth. Recall that $d$ is the dimension of the parameter $\vartheta$.

ASSUMPTION K. The kernel $k$ is a symmetric and twice continuously differentiable density with compact support $[-1, 1]$.

ASSUMPTION B. The bandwidth $b_n$ satisfies $nb_n^4 \to 0$ and $nb_n^{d_*} \to \infty$ with $d_* = (50 + 20d)/(14 + 5d)$.

The requirement on the bandwidth is satisfied by $b_n \sim n^{-\beta}$ for any $\beta$ satisfying $1/4 < \beta < 1/d_*$. Another possibility is $b_n \sim (n \log(n))^{-1/4}$.

In Theorem 2.1 below we describe expansions of $\psi(h, \tilde f)$ and $\psi(h, \tilde f_w)$. For this, we define, for $h \in \mathcal{H}$, a function $\bar h = \bar h_1 + \cdots + \bar h_m$ by

$$\bar h_j(y_j) = \int \cdots \int h(y_1, \ldots, y_m) \prod_{k \ne j} f(y_k) \, dy_k.$$

Note that $\bar h_j(\varepsilon_j) = E(h(\varepsilon_1, \ldots, \varepsilon_m)|\varepsilon_j)$. For a measurable function $g$, we define the $V$-norm by

$$\|g\|_V = \int V(y)|g(y)| \, dy.$$

It follows from Assumptions V and F that $f'$ has finite $V$-norm. Indeed, one has

$$\|f'\|_V^2 = (E[V(\varepsilon)|\ell(\varepsilon)|])^2 \le E[V^2(\varepsilon)] E[\ell^2(\varepsilon)].$$

Recall that

$$\hat{f}(y) = \frac{1}{n}\sum_{i=1}^{n} k_{b_n}(y - \varepsilon_i)$$

denotes the kernel density estimator based on the true innovations. For $g \in L_2(F)$, set

$$B(g) = E[g(\varepsilon)\ell(\varepsilon)]E[\dot{r}_\vartheta(\mathbf{X})],$$
$$U_n(g) = \int g(y)\hat{f}(y)\,dy - \frac{1}{n}\sum_{i=1}^{n} g(\varepsilon_i),$$

let $g^*$ denote the projection of $g$ onto the subspace $\{v \in L_2(F) : \int v(y)f(y)\,dy = 0\}$ and let $g^\#$ denote the projection of $g$ onto the subspace

$$\mathcal{V} = \left\{v \in L_2(F) : \int v(y)f(y)\,dy = \int yv(y)f(y)\,dy = 0\right\}.$$

It is easy to check that $g^*(y) = g(y) - E[g(\varepsilon)]$ and

$$g^\#(y) = g(y) - E[g(\varepsilon)] - \sigma^{-2}E[\varepsilon g(\varepsilon)]\,y, \qquad y \in \mathbb{R}.$$

Since $E[\ell(\varepsilon)] = 0$ and $E[\varepsilon\ell(\varepsilon)] = 1$, we have $\ell^\#(\varepsilon) = \ell(\varepsilon) - \sigma^{-2}\varepsilon$. Note that $E[g(\varepsilon)\ell^\#(\varepsilon)] = E[g^\#(\varepsilon)\ell(\varepsilon)]$. Also, $E[g^*(\varepsilon)\ell(\varepsilon)] = E[g(\varepsilon)\ell(\varepsilon)]$ and $B(g^*) = B(g)$.

Our expansions rely on the following lemma which summarizes results of Müller, Schick and Wefelmeyer [18], namely, their Theorems 3.1–3.3.

LEMMA 2.1. *Suppose Assumptions* B, F, K, R *and* V *hold. Then* $\|\tilde{f} - f\|_V = o_p(n^{-1/4})$ *and* $\|\tilde{f} - \hat{f} - f'E[\dot{r}_\vartheta(\mathbf{X})]^\top(\tilde{\vartheta} - \vartheta)\|_V = o_p(n^{-1/2})$. *Moreover,* $\|\tilde{f}_w - f\|_V = o_p(n^{-1/4})$ *and, with* $\xi(y) = yf(y)$,

$$\left\|\tilde{f}_w - \hat{f} + \sigma^{-2}\xi\frac{1}{n}\sum_{i=1}^{n}\varepsilon_i + \ell^\# fE[\dot{r}_\vartheta(\mathbf{X})]^\top(\tilde{\vartheta} - \vartheta)\right\|_V = o_p(n^{-1/2}).$$

THEOREM 2.1. *Suppose Assumptions* B, F, K, R *and* V *hold. Then*

$$\sup_{h\in\mathcal{H}}\left|\psi(h,\tilde{f}) - \psi(h,f) - \frac{1}{n}\sum_{i=1}^{n}\bar{h}^*(\varepsilon_i) + B(\bar{h}^*)^\top(\tilde{\vartheta} - \vartheta) - U_n(\bar{h})\right|$$
$$= o_p(n^{-1/2});$$
$$\sup_{h\in\mathcal{H}}\left|\psi(h,\tilde{f}_w) - \psi(h,f) - \frac{1}{n}\sum_{i=1}^{n}\bar{h}^\#(\varepsilon_i) + B(\bar{h}^\#)^\top(\tilde{\vartheta} - \vartheta) - U_n(\bar{h})\right|$$
$$= o_p(n^{-1/2}).$$

PROOF. We prove only the second conclusion. For a subset $A$ of $\{1,\dots,m\}$, let

$$\phi_A(\mathbf{y}) = \prod_{j\notin A} f(y_j) \prod_{j\in A} (\tilde f_w(y_j) - f(y_j)), \qquad \mathbf{y} = (y_1,\dots,y_m).$$

Setting $\varphi_r(\mathbf{y}) = \sum_{|A|=r} \phi_A(\mathbf{y})$, we have

$$\prod_{j=1}^m \tilde f_w(y_j) = \prod_{j=1}^m (f(y_j) + \tilde f_w(y_j) - f(y_j)) = \sum_{A\subset\{1,\dots,m\}} \phi_A(\mathbf{y}) = \sum_{r=0}^m \varphi_r(\mathbf{y}).$$

Note that

$$\varphi_0(\mathbf{y}) = \prod_{j=1}^m f(y_j) \quad \text{and} \quad \varphi_1(\mathbf{y}) = \sum_{j=1}^m (\tilde f_w(y_j) - f(y_j)) \prod_{k\neq j}^m f(y_k).$$

Thus,

$$\int h(\mathbf{y})\varphi_0(\mathbf{y})\,d\mathbf{y} = \psi(h,f) \quad \text{and} \quad \int h(\mathbf{y})\varphi_1(\mathbf{y})\,d\mathbf{y} = \int \bar h(y)(\tilde f_w(y) - f(y))\,dy.$$

Using (2.1), we obtain

$$\sum_{r=2}^m \left| \int h(\mathbf{y})\varphi_r(\mathbf{y})\,d\mathbf{y} \right| \le \sum_{r=2}^m \int H(\mathbf{y})|\varphi_r(\mathbf{y})|\,d\mathbf{y} = \sum_{r=2}^m \binom{m}{r} \|\tilde f_w - f\|_V^r \|f\|_V^{m-r}.$$

Since $\|\tilde f_w - f\|_V = o_p(n^{-1/4})$, we obtain

$$\sup_{h\in\mathcal{H}} \left| \psi(h,\tilde f_w) - \psi(h,f) - \int \bar h(y)(\tilde f_w(y) - f(y))\,dy \right| = o_p(n^{-1/2}).$$

Note that $|\bar h| \le C_m V$ with $C_m = m\|f\|_V^{m-1}$. Thus, by the last assertion of Lemma 2.1, $\sup_{h\in\mathcal{H}} |R_n(h)| = o_p(n^{-1/2})$, where

$$R_n(h) = \int \bar h(y) \Bigg( \tilde f_w(y) - \hat f(y) + \sigma^{-2}\xi(y) \frac{1}{n} \sum_{i=1}^n \varepsilon_i + \ell^{\#}(y) f(y) E[\dot r_\vartheta(\mathbf{X})]^\top (\tilde\vartheta - \vartheta) \Bigg)\,dy.$$

Since $E[\bar h(\varepsilon)\ell^{\#}(\varepsilon)]E[\dot r_\vartheta(\mathbf{X})] = B(\bar h^{\#})$, the desired result follows. □

In order to obtain functional central limit theorems for the smoothed von Mises process $\{n^{1/2}(\psi(h,\tilde f) - \psi(h,f)) : h\in\mathcal{H}\}$ based on the residuals and for its weighted version, we can now apply results on smoothed empirical processes $\{n^{1/2}\int g(y)(\hat f(y) - f(y))\,dy : g\in\mathcal{G}\}$ based on the innovations. This also requires an estimator $\tilde\vartheta$ that is *asymptotically linear* in the sense that

$$\tilde\vartheta = \vartheta + \frac{1}{n}\sum_{i=1}^n \varphi(\mathbf{X}_{i-1},\varepsilon_i) + o_p(n^{-1/2}) \tag{2.4}$$

with *influence function* $\varphi(\mathbf{X},\varepsilon)$ satisfying $E(\varphi(\mathbf{X},\varepsilon)|\mathbf{X})=0$ and $E[|\varphi(\mathbf{X},\varepsilon)|^2]<\infty$. Typically, $\varphi$ is orthogonal to $\mathcal{V}$ in the sense that $E[\varphi(\mathbf{X},\varepsilon)v(\varepsilon)]=0$ for all $v\in\mathcal{V}$.

In the literature one decomposes $n^{1/2}\int g(y)(\hat f(y)-f(y))\,dy$ into a variance term

$$n^{1/2}\int g(y)(\hat f(y)-f*k_{b_n}(y))\,dy$$

and a bias term

$$n^{1/2}\int g(y)(f*k_{b_n}(y)-f(y))\,dy.$$

One assumes that the bias term tends to zero uniformly in $g$,

$$\sup_{g\in\mathcal{G}}\left|n^{1/2}\int g(y)(f*k_{b_n}(y)-f(y))\,dy\right|\to 0. \tag{2.5}$$

Sufficient conditions for this analytic property are easily given in terms of smoothness of $f$ and an appropriate bandwidth $b_n$. For example, (2.5) holds if $nb_n^4\to 0$ and

$$\sup_{g\in\mathcal{G}}\left|\int g(y)(f(y-s)-f(y)+sf'(y))\,dy\right|=O(s^2).$$

To deal with the variance term, van der Vaart ([38], (1.1)) and Rost ([25], (2.7)) use a condition that in our case is

$$\sup_{g\in\mathcal{G}}\int\left(\int(g(y+b_nu)-g(y))k(u)\,du\right)^2 f(y)\,dy\to 0. \tag{2.6}$$

van der Vaart [38] shows that if $\mathcal{G}$ is Donsker and translation invariant, then conditions (2.5) and (2.6) imply that the smoothed empirical process converges weakly in $\ell_\infty(\mathcal{G})$ to a tight Brownian bridge process. Inspection of his proof shows that we can remove translation invariance if we strengthen $\mathcal{G}$ being Donsker to $\mathcal{G}_\eta=\{g(\cdot+t):|t|\le\eta, g\in\mathcal{G}\}$ being Donsker for some $\eta>0$.

Suppose now that $\mathcal{G}$ has an envelope $V\in L_2(F)$ satisfying

$$\int(V(y+s)-V(y))^2 f(y)\,dy\to 0 \qquad \text{as } s\to 0. \tag{2.7}$$

Then condition (2.6) holds if $\mathcal{G}$ is totally bounded in $L_2(F)$. This follows from the characterization of compact subsets of $L_2(\nu)$ for finite measures $\nu$ with Lebesgue density given in Lemma 7.1. If $\mathcal{G}$ is Donsker, then $\mathcal{G}$ is totally bounded in $L_2(F)$ and, hence, condition (2.6) holds. We therefore obtain the following version of the Theorem in [38].

PROPOSITION 2.1. *Suppose $\mathcal{G}_\eta$ is Donsker for some $\eta > 0$ and has envelope $V \in L_2(F)$ satisfying condition* (2.7). *Then condition* (2.5) *implies that*

$$\sup_{g\in\mathcal{G}} |U_n(g)| = \sup_{g\in\mathcal{G}} \left| \int g(y)\hat{f}(y)\,dy - \frac{1}{n}\sum_{i=1}^n g(\varepsilon_i) \right| = o_p(n^{-1/2}) \tag{2.8}$$

*and that the smoothed empirical process converges weakly in $\ell_\infty(\mathcal{G})$ to a tight Brownian bridge process.*

One can derive from [25] that $\mathcal{G}_\eta$ Donsker can be replaced by the condition that $\mathcal{G}$ has uniformly integrable $L_2$-entropy. In his Theorem 2.2, Rost uses condition (2.5) with $\mathcal{G}$ replaced by $\mathcal{G} \cup \{V^3\}$ and (2.6). Since $\mathcal{G}$ is totally bounded in $L_2(F)$ if it has uniformly integrable $L_2$-entropy, condition (2.6) is implied by (2.7). Condition (2.5) with $\mathcal{G} = \{V^3\}$ is used only to conclude that $\int V(\varepsilon + b_n u)k(u)\,du$ is uniformly integrable. But the latter follows from condition (2.7). Hence, we have the following version of Rost's Theorem 2.2.

PROPOSITION 2.2. *If $\mathcal{G}$ has uniformly integrable $L_2$-entropy and envelope $V \in L_2(F)$ satisfying* (2.7), *then condition* (2.5) *implies* (2.8) *and the smoothed empirical process converges weakly in $\ell_\infty(\mathcal{G})$ to a tight Brownian bridge process.*

We can now combine Theorem 2.1 and Proposition 2.1 to obtain functional central limit theorems for the von Mises statistics $\psi(h, \tilde{f})$ and $\psi(h, \tilde{f}_w)$. We consider only the weighted version, $\psi(h, \tilde{f}_w)$. Assume that $\tilde{\vartheta}$ is asymptotically linear in the sense of (2.4), with influence function $\varphi$ orthogonal to $\mathcal{V}$. By Theorem 2.1 and Proposition 2.1, $\psi(h, \tilde{f}_w)$ is uniformly asymptotically linear,

$$\sup_{h\in\mathcal{H}} \left| \psi(h, \tilde{f}_w) - \psi(h, f) - \frac{1}{n}\sum_{i=1}^n s_h(\mathbf{X}_{i-1}, \varepsilon_i) \right| = o_p(n^{-1/2}),$$

with influence function $s_h(\mathbf{X}, \varepsilon) = \bar{h}^\#(\varepsilon) - B(\bar{h}^\#)^\top \varphi(\mathbf{X}, \varepsilon)$.

It follows that $\{n^{1/2}(\psi(h, \tilde{f}_w) - \psi(h, f)) : h \in \mathcal{H}\}$ converges weakly in $\ell_\infty(\mathcal{H})$ to a centered Gaussian process with covariance function

$$\begin{aligned} \operatorname{Cov}(h, k) &= E[s_h(\mathbf{X}, \varepsilon)s_k(\mathbf{X}, \varepsilon)] \\ &= E[\bar{h}^\#(\varepsilon)\bar{k}^\#(\varepsilon)] + B(\bar{h}^\#)^\top E[\varphi(\mathbf{X}, \varepsilon)\varphi(\mathbf{X}, \varepsilon)^\top] B(\bar{k}^\#). \end{aligned}$$

We have

$$\begin{aligned} E[\bar{h}^\#(\varepsilon)\bar{k}^\#(\varepsilon)] &= E[\bar{h}(\varepsilon)\bar{k}(\varepsilon)] - E[\bar{h}(\varepsilon)]E[\bar{k}(\varepsilon)] - \sigma^{-2}E[\varepsilon\bar{h}(\varepsilon)]E[\varepsilon\bar{k}(\varepsilon)], \\ B(\bar{h}^\#) &= (E[\bar{h}(\varepsilon)\ell(\varepsilon)] - \sigma^{-2}E[\varepsilon\bar{h}(\varepsilon)])E[\dot{r}_\vartheta(\mathbf{X})]. \end{aligned}$$

A $n^{1/2}$-consistent estimator of the covariance function is obtained using residual-based empirical estimators for $E[\bar{h}^{\#}(\varepsilon)\bar{k}^{\#}(\varepsilon)]$ and $B(\bar{h}^{\#})$ and an appropriate estimator of the asymptotic variance of $\tilde{\vartheta}$. Note that the term of the form $E[h(\varepsilon)\ell(\varepsilon)]$ could be written $\partial_{s=0}E[h(\varepsilon+s)]$, so estimation of $\ell$ could be avoided.

In our applications to estimation of distribution functions and conditional expectations, the class $\mathcal{H}$ consists of functions that may depend on $\vartheta$ and other parameters. To treat the different cases economically, we now formulate a version of Theorem 2.1 for such classes. Suppose that $\mathcal{H}$ is of the form

$$\mathcal{H}^* = \{h_{\tau,q} : |\tau-\vartheta| \le \Delta, q \in \mathcal{Q}\}$$

for some index set $\mathcal{Q}$, and set $\bar{\mathcal{H}}^* = \{\bar{h} : h \in \mathcal{H}^*\}$.

THEOREM 2.2. *Suppose that Assumptions B, F, K and R hold, that $\mathcal{H}^*$ has envelope $H$ of the form $H(y_1,\dots,y_m) = V(y_1)\cdots V(y_m)$ with $V$ satisfying Assumption V, and that*

$$\sup_{q\in\mathcal{Q}}\left|\psi(h_{\tau,q},f) - \psi(h_{\vartheta,q},f) - \Psi_{\vartheta,q}^\top(\tau-\vartheta)\right| = o(|\tau-\vartheta|) \tag{2.9}$$

*for some vector $\Psi_{\vartheta,q}$. Let $\bar{\mathcal{H}}^*_\eta = \{\bar{h}(\cdot+s) : |s| \le \eta, h \in \mathcal{H}^*\}$ be Donsker for some $\eta > 0$;*

$$\sup_{q\in\mathcal{Q}}\int (\bar{h}_{\tau,q}(y) - \bar{h}_{\vartheta,q}(y))^2 f(y)\,dy \to 0 \qquad \text{as } \tau\to\vartheta; \tag{2.10}$$

$$\sup_{|\tau-\vartheta|\le\Delta}\sup_{q\in\mathcal{Q}}\left|\int \bar{h}_{\tau,q}(y)(f(y-s)-f(y)+sf'(y))\,dy\right| = O(s^2). \tag{2.11}$$

*Set $D^*_{\vartheta,q} = \Psi_{\vartheta,q} - B(\bar{h}^*_{\vartheta,q})$ and $D^{\#}_{\vartheta,q} = \Psi_{\vartheta,q} - B(\bar{h}^{\#}_{\vartheta,q})$. Then*

$$\sup_{q\in\mathcal{Q}}\left|\psi(h_{\tilde{\vartheta},q},\tilde{f}) - \psi(h_{\vartheta,q},f) - \frac{1}{n}\sum_{i=1}^n \bar{h}^*_{\vartheta,q}(\varepsilon_i) - (D^*_{\vartheta,q})^\top(\tilde{\vartheta}-\vartheta)\right|$$
$$= o_p(n^{-1/2})$$

*and*

$$\sup_{q\in\mathcal{Q}}\left|\psi(h_{\tilde{\vartheta},q},\tilde{f}_w) - \psi(h_{\vartheta,q},f) - \frac{1}{n}\sum_{i=1}^n \bar{h}^{\#}_{\vartheta,q}(\varepsilon_i) - (D^{\#}_{\vartheta,q})^\top(\tilde{\vartheta}-\vartheta)\right|$$
$$= o_p(n^{-1/2}).$$

*In particular, if $\tilde{\vartheta}$ is asymptotically linear with influence function $\varphi$ orthogonal to $\mathcal{V}$, then the process $\{n^{1/2}(\psi(h_{\tilde{\vartheta},q},\tilde{f}_w) - \psi(h_{\vartheta,q},f)) : q \in \mathcal{Q}\}$ converges weakly in $\ell_\infty(\mathcal{Q})$ to a centered Gaussian process with covariance function*

$$\mathrm{Cov}(p,q) = E[\bar{h}^{\#}_{\vartheta,p}(\varepsilon)\bar{h}^{\#}_{\vartheta,q}(\varepsilon)] + (D^{\#}_{\vartheta,p})^\top E[\varphi(\mathbf{X},\varepsilon)\varphi(\mathbf{X},\varepsilon)^\top]D^{\#}_{\vartheta,q}.$$

Proof. We prove only the second expansion. It follows from (2.9) and the $n^{1/2}$-consistency of $\tilde{\vartheta}$ that

$$\sup_{q\in\mathcal{Q}} |\psi(h_{\tilde{\vartheta},q}, f) - \psi(h_{\vartheta,q}, f) - \Psi_{\vartheta,q}^\top(\tilde{\vartheta} - \vartheta)| = o_p(n^{-1/2}).$$

It follows from (2.10) that

$$\sup_{q\in\mathcal{Q}} |B(\bar{h}^{\#}_{\tilde{\vartheta},q}) - B(\bar{h}^{\#}_{\vartheta,q})| = o_p(1).$$

Since $\bar{\mathcal{H}}^*$ is Donsker, we obtain from (2.10) that

$$\sup_{q\in\mathcal{Q}} \left| \int (\bar{h}_{\tilde{\vartheta},q}(y) - \bar{h}_{\vartheta,q}(y))\, d(\hat{F}(y) - F(y)) \right| = o_p(n^{-1/2}),$$

with $\hat{F}(y) = \frac{1}{n}\sum_{i=1}^n \mathbf{1}[\varepsilon_i \le y]$. Since $\bar{\mathcal{H}}^*_\eta$ is Donsker and (2.11) holds, we obtain from Proposition 2.1 that

$$\sup_{q\in\mathcal{Q}} \left| \int \bar{h}_{\tilde{\vartheta},q}(y)\hat{f}(y)\, dy - \int \bar{h}_{\tilde{\vartheta},q}(y)\, d\hat{F}(y) \right| = o_p(n^{-1/2}).$$

The desired result now follows from Theorem 2.1. □

A sufficient condition for (2.11) is $\|f(\cdot - s) - f + sf'\|_V = O(s^2)$. This holds, for example, if $\|f'(\cdot - s) - f'\|_V = O(s)$. In particular, it holds if $f'$ is absolutely continuous with $\|f''\|_V$ finite.

Also of interest is the case when $\mathcal{H} = \{h_q : q \in \mathcal{Q}\}$. In this case, the assumptions of Theorem 2.2 simplify considerably.

Corollary 2.1. *Suppose that Assumptions B, F, K and R hold and $\mathcal{H} = \{h_q : q \in \mathcal{Q}\}$ has envelope $H$ of the form $H(y_1, \dots, y_m) = V(y_1)\cdots V(y_m)$ with $V$ satisfying Assumption V. Let $\bar{\mathcal{H}}_\eta = \{\bar{h}_q(\cdot + s) : |s| \le \eta, q \in \mathcal{Q}\}$ be Donsker for some $\eta > 0$ and*

$$\sup_{q\in\mathcal{Q}} \left| \int \bar{h}_q(y)(f(y - s) - f(y) + sf'(y))\, dy \right| = O(s^2). \tag{2.12}$$

*Then*

$$\sup_{q\in\mathcal{Q}} \left| \psi(h_q, \tilde{f}) - \psi(h_q, f) - \frac{1}{n}\sum_{i=1}^n \bar{h}^*_q(\varepsilon_i) - B(\bar{h}^*_q)^\top(\tilde{\vartheta} - \vartheta) \right| = o_p(n^{-1/2});$$

$$\sup_{q\in\mathcal{Q}} \left| \psi(h_q, \tilde{f}_w) - \psi(h_q, f) - \frac{1}{n}\sum_{i=1}^n \bar{h}^{\#}_q(\varepsilon_i) - B(\bar{h}^{\#}_q)^\top(\tilde{\vartheta} - \vartheta) \right| = o_p(n^{-1/2}).$$

*In particular, if $\tilde{\vartheta}$ is asymptotically linear with influence function $\varphi$ orthogonal to $\mathcal{V}$, then the process $\{n^{1/2}(\psi(h_q, \tilde{f}_w) - \psi(h_q, f)) : q \in \mathcal{Q}\}$ converges weakly in $\ell_\infty(\mathcal{Q})$ to a centered Gaussian process with covariance function*

$$\mathrm{Cov}(p, q) = E[\bar{h}^{\#}_p(\varepsilon)\bar{h}^{\#}_q(\varepsilon)] + B(\bar{h}^{\#}_p)^\top E[\varphi(\mathbf{X}, \varepsilon)\varphi(\mathbf{X}, \varepsilon)^\top] B(\bar{h}^{\#}_q).$$

**3. Conditional expectations of lag one.** Let $\mathcal{Q}$ be a family of functions from $\mathbb{R}$ to $\mathbb{R}$. For $q \in \mathcal{Q}$, the conditional expectation $E(q(X_{n+1})|\mathbf{X}_n = \mathbf{x})$ can be written as $\nu(\vartheta, q) = E[q(\varepsilon + r_\vartheta(\mathbf{x}))]$. We estimate $\nu(\vartheta, q)$ by

$$\tilde{\nu}(q) = \int q(y + r_{\tilde{\vartheta}}(\mathbf{x}))\tilde{f}(y)\,dy$$

and

$$\tilde{\nu}_w(q) = \int q(y + r_{\tilde{\vartheta}}(\mathbf{x}))\tilde{f}_w(y)\,dy.$$

THEOREM 3.1. *Suppose* $\mathcal{G}_\eta = \{q(\cdot + r_\vartheta(\mathbf{x}) + s) : |s| \le \eta, q \in \mathcal{Q}\}$ *is Donsker for some* $\eta > 0$ *and has an envelope* $V$ *that satisfies Assumption* V. *Suppose* $f$ *has finite Fisher information for location and satisfies*

$$\sup_{|t|\le\eta}\sup_{q\in\mathcal{Q}}\left|\int q(y + r_\vartheta(\mathbf{x}) + t)(f(y-s) - f(y) + sf'(y))\,dy\right| = O(s^2). \tag{3.1}$$

*Let Assumptions* B, K *and* R *hold. Then*

$$\sup_{q\in\mathcal{Q}}\left|\tilde{\nu}(q) - \frac{1}{n}\sum_{i=1}^n q(\varepsilon_i + r_\vartheta(\mathbf{x})) - D_q^\top(\tilde{\vartheta} - \vartheta)\right| = o_p(n^{-1/2});$$

$$\sup_{q\in\mathcal{Q}}\left|\tilde{\nu}_w(q) - \frac{1}{n}\sum_{i=1}^n (q(\varepsilon_i + r_\vartheta(\mathbf{x})) - c_q\varepsilon_i) - \bar{D}_q^\top(\tilde{\vartheta} - \vartheta)\right| = o_p(n^{-1/2}),$$

*where* $D_q = E[q(\varepsilon + r_\vartheta(\mathbf{x}))\ell(\varepsilon)](\dot{r}_\vartheta(\mathbf{x}) - E[\dot{r}_\vartheta(\mathbf{X})])$ *and* $\bar{D}_q = D_q + c_q E[\dot{r}_\vartheta(\mathbf{X})]$ *with* $c_q = \sigma^{-2}E[\varepsilon q(\varepsilon + r_\vartheta(\mathbf{x}))]$.

*In particular, if* $\tilde{\vartheta}$ *is asymptotically linear with influence function* $\varphi$ *orthogonal to* $\mathcal{V}$, *then the process* $\{n^{1/2}(\tilde{\nu}_w(q) - E[q(\varepsilon + r_\vartheta(\mathbf{x}))]) : q \in \mathcal{Q}\}$ *converges weakly in* $\ell_\infty(\mathcal{Q})$ *to a centered Gaussian process with covariance function*

$$\begin{aligned}\mathrm{Cov}(p, q) = {} & E[p(\varepsilon + r_\vartheta(\mathbf{x}))q(\varepsilon + r_\vartheta(\mathbf{x}))] - E[p(\varepsilon + r_\vartheta(\mathbf{x}))]E[q(\varepsilon + r_\vartheta(\mathbf{x}))] \\ & - \sigma^2 c_p c_q + D_p^\top E[\varphi(\mathbf{X}, \varepsilon)\varphi(\mathbf{X}, \varepsilon)^\top]D_q.\end{aligned}$$

PROOF. We apply Theorem 2.2 with

$$\bar{\mathcal{H}}^* = \mathcal{H}^* = \{q(\cdot + r_\tau(\mathbf{x})) : |\tau - \vartheta| \le \Delta, q \in \mathcal{Q}\}$$

and some small positive $\Delta$. In view of Assumption R, we can take $\Delta$ sufficiently small for $\mathcal{H}^*_{\eta/2}$ to be contained in $\mathcal{G}_\eta$. Thus, condition (2.11) is implied by (3.1). Since

$$\begin{aligned}&\int (q(y + r_\tau(\mathbf{x})) - q(y + r_\vartheta(\mathbf{x})))f(y)\,dy \\ &\qquad = \int q(y + r_\vartheta(\mathbf{x}))(f(y - (r_\tau(\mathbf{x}) - r_\vartheta(\mathbf{x}))) - f(y))\,dy,\end{aligned}$$

it follows from (3.1), differentiability of $\tau \mapsto r_\tau(\mathbf{x})$ at $\vartheta$ and finiteness of $\|f'\|_V$ that condition (2.9) holds for the present $\mathcal{H}^*$ with

$$\Psi_{\vartheta,q} = E[q(\varepsilon + r_\vartheta(\mathbf{x}))\ell(\varepsilon)]\dot{r}_\vartheta(\mathbf{x}).$$

As $\mathcal{G} = \{q(\cdot + r_\vartheta(\mathbf{x})) : q \in \mathcal{Q}\}$ is totally bounded and has an envelope $V$ that satisfies (2.7), condition (2.10) is met by the compactness criterion in Lemma 7.1. □

The conditional distribution function of $X_{n+1}$, given $\mathbf{X}_n = \mathbf{x}$, can be written $t \mapsto F(t - r_\vartheta(\mathbf{x}))$. We can estimate it by $\tilde{F}_s(t - r_{\tilde{\vartheta}}(\mathbf{x}))$ or $\tilde{F}_{sw}(t - r_{\tilde{\vartheta}}(\mathbf{x}))$, where $\tilde{F}_s$ and $\tilde{F}_{sw}$ are the distribution functions corresponding to $\tilde{f}$ and $\tilde{f}_w$, respectively. The corresponding class $\mathcal{Q}$ is $\{\mathbf{1}_{(-\infty,t]} : t \in \mathbb{R}\}$; it is Donsker and translation invariant. Its envelope is $V = 1$, which satisfies Assumption V. Here the left-hand side of (3.1) becomes $\sup_{s\in\mathbb{R}} |F(t-s) - F(t) + sf(t)|$. Thus, (3.1) holds if $f$ is Lipschitz. Hence, Theorem 3.1 implies the following result.

THEOREM 3.2. *Suppose Assumptions B, K and R hold. Let $f$ be Lipschitz and have finite Fisher information for location. Then*

$$\sup_{t\in\mathbb{R}} \left| \tilde{F}_s(t - r_{\tilde{\vartheta}}(\mathbf{x})) - \frac{1}{n}\sum_{i=1}^n \mathbf{1}[\varepsilon_i \le t - r_\vartheta(\mathbf{x})] - D_t^\top(\tilde{\vartheta} - \vartheta) \right|$$
$$= o_p(n^{-1/2});$$
$$\sup_{t\in\mathbb{R}} \left| \tilde{F}_{sw}(t - r_{\tilde{\vartheta}}(\mathbf{x})) - \frac{1}{n}\sum_{i=1}^n (\mathbf{1}[\varepsilon_i \le t - r_\vartheta(\mathbf{x})] - c_t\varepsilon_i) - \bar{D}_t^\top(\tilde{\vartheta} - \vartheta) \right|$$
$$= o_p(n^{-1/2}),$$

*where $D_t = -f(t - r_\vartheta(\mathbf{x}))(\dot{r}_\vartheta(\mathbf{x}) - E[\dot{r}_\vartheta(\mathbf{X})])$ and $\bar{D}_t = D_t + c_t E[\dot{r}_\vartheta(\mathbf{X})]$ with*

$$c_t = \sigma^{-2} \int_{-\infty}^{t - r_\vartheta(\mathbf{x})} y f(y)\, dy.$$

*In particular, if $\tilde{\vartheta}$ is asymptotically linear with influence function $\varphi$ orthogonal to $\mathcal{V}$, then the process $\{n^{1/2}(\tilde{F}_{sw}(t - r_{\tilde{\vartheta}}(\mathbf{x})) - F(t - r_\vartheta(\mathbf{x}))) : t \in \mathbb{R}\}$ converges weakly in $\ell_\infty(\mathbb{R})$ to a centered Gaussian process with covariance function*

$$\operatorname{Cov}(s,t) = F((s - r_\vartheta(\mathbf{x})) \wedge (t - r_\vartheta(\mathbf{x}))) - F(s - r_\vartheta(\mathbf{x}))F(t - r_\vartheta(\mathbf{x}))$$
$$- \sigma^2 c_s c_t + \bar{D}_s^\top E[\varphi(\mathbf{X},\varepsilon)\varphi(\mathbf{X},\varepsilon)^\top]\bar{D}_t.$$

EXAMPLE 3.1. For $0 < u < 1$, let $\psi_u(G) = G^{-1}(u) = \inf\{t : G(t) \ge u\}$ denote the left-inverse of a distribution function $G$ at $u$. The *conditional* *u-quantile* of $X_{n+1}$, given $\mathbf{X}_n = \mathbf{x}$, is $\psi_u(F(\cdot - r_\vartheta(\mathbf{x})) = F^{-1}(u) + r_\vartheta(\mathbf{x})$. We can estimate it by $\psi_u(\tilde{F}_{sw}(\cdot - r_{\tilde\vartheta}(\mathbf{x}))) = \tilde{F}_{sw}^{-1}(u) + r_{\tilde\vartheta}(\mathbf{x})$. Let $0 < c \le d < 1$. Recall that we assumed that the density $f$ is positive. Thus, by Proposition 1 of [10] on compact differentiability of quantile functions, we obtain the uniform stochastic expansion

$$\sup_{u\in[c,d]}\Bigg| \tilde{F}_{sw}^{-1}(u) + r_{\tilde\vartheta}(\mathbf{x}) - (F^{-1}(u) + r_\vartheta(\mathbf{x})) + \frac{1}{f(F^{-1}(u))}\frac{1}{n}\sum_{i=1}^n(\mathbf{1}[\varepsilon_i \le F^{-1}(u)] - u - a_u\varepsilon_i) + \left(\left(1+\frac{a_u}{f(F^{-1}(u))}\right)E[\dot r_\vartheta(\mathbf{X})] - \dot r_\vartheta(x)\right)^\top(\tilde\vartheta - \vartheta)\Bigg| = o_p(n^{-1/2})$$

with

$$a_u = \sigma^{-2}\int_{-\infty}^{F^{-1}(u)} yf(y)\,dy.$$

It follows that the smoothed and weighted conditional quantile process

$$\{n^{1/2}(\tilde{F}_{sw}^{-1}(u) + r_{\tilde\vartheta}(\mathbf{x}) - (F^{-1}(u) + r_\vartheta(\mathbf{x}))) : u \in [c,d]\}$$

converges weakly in $\ell_\infty([c,d])$ to a centered Gaussian process.

EXAMPLE 3.2. Consider the classical AR(1) model $X_i = \vartheta X_{i-1} + \varepsilon_i$ with $|\vartheta| < 1$. It satisfies Condition R with $\dot r_\vartheta(x) = x$ and $E[\dot r_\vartheta(X)] = 0$. A natural estimator for $\vartheta$ is the least squares estimator $\tilde\vartheta$, which has expansion

$$\tilde\vartheta = \frac{\sum_{i=1}^n X_{i-1}X_i}{\sum_{i=1}^n X_{i-1}^2} = \vartheta + \frac{1}{n}\sum_{i=1}^n \frac{1-\vartheta^2}{\sigma^2}X_{i-1}\varepsilon_i + o_p(n^{-1/2}). \tag{3.2}$$

Fix $t$ and $x$ in $\mathbb{R}$. For the estimator of the conditional distribution function at $t$ of $X_{n+1}$, given $X_n = x$, we obtain

$$\tilde{F}_s(t - \tilde\vartheta x) = \frac{1}{n}\sum_{i=1}^n\left(\mathbf{1}[\varepsilon_i \le t - \vartheta x] - xf(t-\vartheta x)\frac{1-\vartheta^2}{\sigma^2}X_{i-1}\varepsilon_i\right) + o_p(n^{-1/2});$$

$$\tilde{F}_{sw}(t - \tilde\vartheta x) = \frac{1}{n}\sum_{i=1}^n\left(\mathbf{1}[\varepsilon_i \le t - \vartheta x] - c_t\varepsilon_i - xf(t-\vartheta x)\frac{1-\vartheta^2}{\sigma^2}X_{i-1}\varepsilon_i\right) + o_p(n^{-1/2}).$$

It follows that $n^{1/2}(\tilde{F}_s(t - \tilde\vartheta x) - F(t - \vartheta x))$ is asymptotically normal with mean zero and variance $\tau^2 = F(t-\vartheta x)(1 - F(t-\vartheta x)) + x^2 f^2(t-\vartheta x)(1-\vartheta^2)$, while $n^{1/2}(\tilde{F}_{sw}(t-\tilde{\vartheta}x) - F(t-\vartheta x))$ is asymptotically normal with mean zero and variance $\tau^2 - \sigma^2 c_t^2$. Thus, weighting results in a smaller asymptotic variance. For $t = x = 0$ and $f$ the standard normal density, the asymptotic variances are $1/4$ and $1/4 - 1/(2\pi) \simeq 0.0908$. In this case weighting reduces the asymptotic variance by about 64%.

Now consider the case where $\mathcal{Q}$ consists of one element $q$. The corresponding class $\mathcal{G}_\eta$ equals $\{q(\cdot + r_\vartheta(\mathbf{x}) + s) : |s| \le \eta\}$. Assume now that $f$ has a finite absolute moment of order greater than $2\gamma + 1$ and that $q$ satisfies the growth condition

$$|q(y)| \le (1 + |y|)^\gamma, \qquad y \in \mathbb{R}, \tag{3.3}$$

and the Lipschitz condition

$$|q(y + s_1) - q(y + s_2)| \le L|s_1 - s_2|(1 + |y|)^\gamma, \qquad y \in \mathbb{R}, \tag{3.4}$$

for $s_1$, $s_2$ in a neighborhood of $r_\vartheta(\mathbf{x})$. Then $\mathcal{G}_\eta$ has envelope $V$ of the form $V(y) = K(1 + |y|)^\gamma$, which satisfies Assumption V. Also, $\mathcal{G}_\eta$ is Donsker. This follows since the bracketing numbers $N_{[\cdot]}(\delta, \mathcal{G}_\eta, L_2(F))$ are of order $1/\delta$; take brackets of the form $q(\cdot + r_\vartheta(\mathbf{x}) + s_j) \mp c\delta V$. The left-hand side of (3.1) now becomes $\sup_{|t|<\eta} |\Delta_{s,t}|$ with

$$\Delta_{s,t} = \int q(y + r_\vartheta(\mathbf{x}) + t)(f(y - s) - f(y) + sf'(y))\,dy.$$

We can write

$$\begin{aligned}
\Delta_{s,t} &= -\int q(y + r_\vartheta(\mathbf{x}) + t) \int_0^1 s(f'(y - us) - f'(y))\,du\,dy \\
&= -s \int_0^1 \int (q(y + r_\vartheta(\mathbf{x}) + t + us) - q(y + r_\vartheta(\mathbf{x}) + t)) f'(y)\,dy\,du.
\end{aligned}$$

By the Lipschitz property of $q$ and the finiteness of $\|f'\|_V$ under Assumption F, we obtain $\sup_{|t|<\eta} |\Delta_{s,t}| = O(s^2)$, which is (3.1). Thus, Theorem 3.1 implies the following result.

THEOREM 3.3. *Suppose Assumptions* B, K *and* R *hold, $q$ satisfies* (3.3) *and* (3.4) *and $f$ has finite Fisher information for location and finite absolute moment of order greater than $2\gamma + 1$. Then*

$$\int q(y + r_{\tilde{\vartheta}}(\mathbf{x}))\tilde{f}(y)\,dy = \frac{1}{n}\sum_{i=1}^n q(\varepsilon_i + r_\vartheta(\mathbf{x})) + D_q^\top(\tilde{\vartheta} - \vartheta) + o_p(n^{-1/2});$$

$$\begin{aligned}
\int q(y + r_{\tilde{\vartheta}}(\mathbf{x}))\tilde{f}_w(y)\,dy = \frac{1}{n}\sum_{i=1}^n (q(\varepsilon_i + r_\vartheta(\mathbf{x})) - c_q\varepsilon_i) + \bar{D}_q^\top(\tilde{\vartheta} - \vartheta) \\
+ o_p(n^{-1/2}),
\end{aligned}$$

*with $c_q$, $D_q$ and $\bar{D}_q$ as in Theorem* 3.1.

Theorem 3.3 can be used to estimate conditional moments and absolute moments of lag one. For example, to treat estimation of the conditional $\gamma$th absolute moment $E(|X_{n+1}|^\gamma|\mathbf{X}_n = \mathbf{x})$ with $\gamma \geq 1$, take $q(y) = |y|^\gamma$. Our estimators are $\int |y + r_{\tilde{\vartheta}}(\mathbf{x})|^\gamma \tilde{f}(y)\,dy$ and its weighted version $\int |y + r_{\tilde{\vartheta}}(\mathbf{x})|^\gamma \tilde{f}_w(y)\,dy$.

**4. Conditional expectations of lag two.** Let $\mathcal{Q}$ be a family of functions from $\mathbb{R}^2$ to $\mathbb{R}$. For $q \in \mathcal{Q}$, the conditional expectation

$$E(q(X_{n+1}, X_{n+2})|\mathbf{X}_n = \mathbf{x})$$

can be written

$$\nu(\vartheta, q) = E[q(\varrho_\vartheta(\varepsilon_1, \varepsilon_2))] = \iint q(\varrho_\vartheta(y, z)) f(y) f(z)\,dy\,dz$$

with

$$\varrho_\vartheta(y, z) = (y + r_\vartheta(\mathbf{x}), z + r_\vartheta(\mathbf{x}_{-1}, y + r_\vartheta(\mathbf{x}))),$$

where $\mathbf{x}_{-1} = (x_2, \ldots, x_p)$. We estimate $\nu(\vartheta, q)$ by

$$\tilde{\nu}(q) = \iint q(\varrho_{\tilde{\vartheta}}(y, z)) \tilde{f}(y) \tilde{f}(z)\,dy\,dz$$

and its weighted version

$$\tilde{\nu}_w(q) = \iint q(\varrho_{\tilde{\vartheta}}(y, z)) \tilde{f}_w(y) \tilde{f}_w(z)\,dy\,dz.$$

We shall apply Theorem 2.2 to obtain stochastic expansions for these estimators.

We have

$$h_{\tau,q}(y, z) = q(\varrho_\tau(y, z)) \qquad \text{and} \qquad \bar{h}_{\tau,q} = \bar{h}^{(1)}_{\tau,q} + \bar{h}^{(2)}_{\tau,q}$$

with

$$\bar{h}^{(1)}_{\tau,q}(y) = \int q(\varrho_\tau(y, u)) f(u)\,du \qquad \text{and} \qquad \bar{h}^{(2)}_{\tau,q}(z) = \int q(\varrho_\tau(u, z)) f(u)\,du.$$

To get an envelope for the class $\mathcal{H}^* = \{h_{\tau,q} : |\tau - \vartheta| \leq \Delta, q \in \mathcal{Q}\}$, we assume that $\mathcal{Q}$ has an envelope $V_\mathcal{Q}$ of the form

$$V_\mathcal{Q}(x_1, x_2) = C_\mathcal{Q}(1 + |x_1|)^{\gamma_1}(1 + |x_2|)^{\gamma_2} \tag{4.1}$$

for some finite constant $C_\mathcal{Q}$ and nonnegative exponents $\gamma_1$ and $\gamma_2$, and impose the following growth condition on the autoregression functions: for some constant $A$,

$$|r_\tau(u_1, \ldots, u_p)| \leq A\left(1 + \sum_{j=1}^p |u_j|\right), \qquad |\tau - \vartheta| \leq \Delta. \tag{4.2}$$

Such a growth condition is typically needed for ergodicity of the model; see [4, 5] and [1]. There is then a constant $C'_{\mathcal{Q}}$ such that

$$
\begin{aligned}
|q(\varrho_\tau(y,z))| &\le C_{\mathcal{Q}}(1+|y+r_\tau(\mathbf{x})|)^{\gamma_1}(1+|z+r_\tau(\mathbf{x}_{-1},y+r_\tau(\mathbf{x}))|)^{\gamma_2}\\
&\le C'_{\mathcal{Q}}(1+|y|)^{\gamma_1}(1+|z|+|y|)^{\gamma_2} \qquad (4.3)\\
&\le C'_{\mathcal{Q}}(1+|y|)^{\gamma_1+\gamma_2}(1+|z|)^{\gamma_2}.
\end{aligned}
$$

Thus, $\mathcal{H}^*$ has an envelope $H$ of the form $H(y,z)=V(y)V(z)$ with $V(y)=K(1+|y|)^{\gamma_1+\gamma_2}$. This $V$ satisfies Assumption V if $f$ has finite absolute moment of order greater than $2\gamma_1+2\gamma_2+1$. We can now use the special structure of $h_{\tau,q}$ to show that (2.9) holds.

Lemma 4.1. *Let $\mathcal{Q}$ have envelope $V_{\mathcal{Q}}$ of the form* (4.1). *Suppose that $f$ has finite Fisher information for location and finite absolute moment of order greater than $2\gamma_1+2\gamma_2+1$. Suppose Assumption* R *and the growth condition* (4.2) *hold and that*

$$
\int (r_{\vartheta+t}(\mathbf{x}_{-1},y)-r_\vartheta(\mathbf{x}_{-1},y)-\dot r_\vartheta(\mathbf{x}_{-1},y)^\top t)^2 f(y-r_\vartheta(\mathbf{x}))\,dy = o(|t|^2). \tag{4.4}
$$

*Then*

$$
\sup_{q\in\mathcal{Q}}\left|\iint (q(\varrho_{\vartheta+t}(y,z))-q(\varrho_\vartheta(y,z))-q(\varrho_\vartheta(y,z))\chi(y,z)^\top t)f(y)f(z)\,dy\,dz\right|
$$
$$
= o(|t|),
$$

*where $\chi(y,z)=\ell(y)\dot r_\vartheta(\mathbf{x})+\ell(z)\dot r_\vartheta(\mathbf{x}_{-1},y+r_\vartheta(\mathbf{x}))$. Thus, condition* (2.9) *holds with*

$$
\begin{aligned}
\Psi_{\vartheta,q} &= \iint q(\varrho_\vartheta(y,z))\chi(y,z)f(y)f(z)\,dy\,dz\\
&= -\iint q(\varrho_\vartheta(y,z))\\
&\qquad \times (f'(y)f(z)\dot r_\vartheta(\mathbf{x})+f(y)f'(z)\dot r_\vartheta(\mathbf{x}_{-1},y+r_\vartheta(\mathbf{x})))\,dy\,dz.
\end{aligned}
$$

Proof. It is easy to check that $\varrho_\tau(\varepsilon_1,\varepsilon_2)$ has a density $p_\tau$ with respect to the Lebesgue measure $\lambda_2$ on $\mathbb{R}^2$ of the form

$$
p_\tau(y,z)=f(y-r_\tau(\mathbf{x}))f(z-r_\tau(\mathbf{x}_{-1},y)). \tag{4.5}
$$

We can write the integral in the assertion as

$$
\iint q(y,z)(p_{\vartheta+t}(y,z)-p_\vartheta(y,z)-\tilde\chi(y,z)^\top t p_\vartheta(y,z))\,dy\,dz, \tag{4.6}
$$

where $\tilde{\chi}$ is the score function at $\vartheta$ of the parametric model $\mathcal{P} = \{p_\tau : |\tau - \vartheta| \le \Delta\}$:

$$\tilde{\chi}(y,z) = \ell(y - r_\vartheta(\mathbf{x}))\dot{r}_\vartheta(\mathbf{x}) + \ell(z - r_\vartheta(\mathbf{x}_{-1}, y))\dot{r}_\vartheta(\mathbf{x}_{-1}, y). \tag{4.7}$$

Actually, $\tilde{\chi}$ is the Hellinger derivative of this model at $\vartheta$. Indeed, since $f$ has finite Fisher information for location, the model $\{f(\cdot - r_\tau(\mathbf{x})) : |\tau - \vartheta| \le \Delta\}$ is Hellinger differentiable at $\vartheta$ with Hellinger derivative $\dot{r}_\vartheta(\mathbf{x})\ell(\cdot - r_\vartheta(\mathbf{x}))$ and this and (4.4) yield the Hellinger differentiability of $\mathcal{P}$ at $\vartheta$ with Hellinger derivative $\tilde{\chi}$; see Proposition A.6 in [13]. It is easy to check that $\int V_{\mathcal{Q}}^2 p_\tau \, d\lambda_2 \to \int V_{\mathcal{Q}}^2 p_\vartheta \, d\lambda_2$ as $\tau \to \vartheta$. Thus, Lemma 7.2 yields the desired result. □

We now address sufficient conditions for (2.10).

Lemma 4.2. *Suppose the assumptions of Lemma 4.1 hold. Then (2.10) is implied by*

$$\sup_{q \in \mathcal{Q}} \int (\bar{h}^{(1)}_{\vartheta,q}(y+s) - \bar{h}^{(1)}_{\vartheta,q}(y))^2 f(y)\, dy \to 0 \qquad \text{as } s \to 0, \tag{4.8}$$

$$\sup_{q \in \mathcal{Q}} \int \left( \int (q(\varrho_\vartheta(y, z + \Delta_\tau(y))) - q(\varrho_\vartheta(y,z))) f(y)\, dy \right)^2 f(z)\, dz \to 0 \qquad \text{as } \tau \to \vartheta, \tag{4.9}$$

*where* $\Delta_\tau(y) = r_\tau(\mathbf{x}_{-1}, y + r_\vartheta(\mathbf{x})) - r_\vartheta(\mathbf{x}_{-1}, y + r_\vartheta(\mathbf{x}))$.

Proof. With $s = r_\tau(\mathbf{x}) - r_\vartheta(\mathbf{x})$, we can write $\varrho_\tau(y,z) = \varrho_\vartheta(y + s, z + \Delta_\tau(y+s))$ and then

$$\bar{h}^{(1)}_{\tau,q}(y) = \int q(\varrho_\vartheta(y+s, z)) f(z - \Delta_\tau(y+s))\, dz,$$

$$\bar{h}^{(2)}_{\tau,q}(y) = \int q(\varrho_\vartheta(y, z + \Delta_\tau(y))) f(y - s)\, dy.$$

In view of (4.8) and (4.9), it suffices to show that, as $\tau \to \vartheta$,

$$\sup_{q \in \mathcal{Q}} \int (\bar{h}^{(1)}_{\tau,q}(y) - \bar{h}^{(1)}_{\vartheta,q}(y+s))^2 f(y)\, dy \to 0, \tag{4.10}$$

$$\sup_{q \in \mathcal{Q}} \int \left( \int q(\varrho_\vartheta(y, z + \Delta_\tau(y)))(f(y-s) - f(y))\, dy \right)^2 f(z)\, dz \to 0. \tag{4.11}$$

The Cauchy–Schwarz inequality gives

$$(\bar{h}^{(1)}_{\tau,q}(y) - \bar{h}^{(1)}_{\vartheta,q}(y+s))^2 = \left( \int q(\varrho_\vartheta(y+s, z))(f(z - \Delta_\tau(y+s)) - f(z))\, dz \right)^2$$

$$\leq \int q^2(\varrho_\vartheta(y+s,z))(f(z-\Delta_\tau(y+s))+f(z))\,dz$$
$$\times \int |f(z-\Delta_\tau(y+s))-f(z)|\,dz.$$

Using (4.3), we can bound the first integral on the right-hand side by $C_{\mathcal{Q}}'^2[(1+|y|)^{2\gamma_1+2\gamma_2}+(1+|y+s|)^{2\gamma_1+2\gamma_2}]\int(1+|z|)^{2\gamma_2}f(z)\,dz$, while the second integral can be bounded by $|\Delta_\tau(y+s)|\|f'\|_1$. Indeed, since $f$ has finite Fisher information, $f'$ is integrable and $\int|f(z-v)-f(z)|\,dz \leq |v|\|f'\|_1$ for every real $v$. Using these bounds, we obtain that the left-hand side of (4.10) is bounded for $|s|<1$ by a constant times

$$\int(1+|y|)^{2\gamma_1+2\gamma_2}|\Delta_\tau(y)|f(y-s)\,dy.$$

Since $|\Delta_\tau(y)| \leq \tilde{A}(1+|y|)$ and $\Delta_\tau(y)\to 0$ for every $y$, and since

$$\int(1+|y|)^{2\gamma_1+2\gamma_2+1}|f(y-s)-f(y)|\,dy \to 0,$$

we get (4.10). A similar argument yields (4.11). □

REMARK 4.1. In view of the characterization of compact subsets of $L_2(F)$ given in Lemma 7.1, the above assumptions imply that condition (4.8) is equivalent to total boundedness of $\mathcal{H}^{(1)}=\{\bar{h}^{(1)}_{\vartheta,q}: q\in\mathcal{Q}\}$ in $L_2(F)$. Consequently, (4.8) holds if $\mathcal{Q}$ is a finite set or if $\mathcal{H}^{(1)}$ is Donsker.

Let us now assume that the class $\mathcal{Q}$ satisfies the Lipschitz property

$$|q(y_1,z_1)-q(y_2,z_2)| \leq L_1(y,z)|y_1-y_2| + L_2(y,z)|z_1-z_2|, \tag{4.12}$$

where $y=|y_1|\vee|y_2|$ and $z=|z_1|\vee|z_2|$ and where

$$L_1(y,z)=C_1(1+|y|)^{\alpha_1}(1+|z|)^{\alpha_2} \quad\text{and}\quad L_2(y,z)=C_2(1+|y|)^{\beta_1}(1+|z|)^{\beta_2}$$

for constants $C_1$, $C_2$ and nonnegative exponents $\alpha_1$, $\alpha_2$, $\beta_1$ and $\beta_2$. Let us set $\zeta=\max\{\alpha_1+\alpha_2,\beta_1+\beta_2\}$. Then we derive that, for each $C$, there is a $C_*$ such that, for all $y$, $z$, all $|s_1|,|s_2|,|t_1|,|t_2|\leq C$ and $|a_1(y)|,|a_2(y)|\leq C(1+|y|)$,

$$|q(y+s_1,z+t_1+a_1(y)) - q(y+s_2,z+t_2+a_2(y))|$$
$$\leq C_*(|s_1-s_2|+|t_1-t_2|+|a_1(y)-a_2(y)|)(1+|y|)^\zeta(1+|z|)^\zeta.$$

With the help of this inequality, it is now easy to check that, under the assumptions of Lemma 4.1, the statements (4.8) and (4.9) are met, so that (2.10) holds by Lemma 4.2. Using

$$f(y-s)-f(y)+sf'(y) = -s\int_0^1 (f'(y-ws)-f'(y))\,dw, \tag{4.13}$$

the left-hand side of (2.11) can be bounded by $|s|(T_1(s)+T_2(s))$, where

$$\begin{aligned} T_1(s) &= \sup_{0\le w\le 1}\sup_{|\tau-\vartheta|\le\Delta}\sup_{q\in\mathcal{Q}}\Big|\int\!\!\int (q(\varrho_\tau(y+ws,z)) - q(\varrho_\tau(y,z))) \\ &\qquad\qquad \times f'(y)f(z)\,dy\,dz\Big|, \\ T_2(s) &= \sup_{0\le w\le 1}\sup_{|\tau-\vartheta|\le\Delta}\sup_{q\in\mathcal{Q}}\Big|\int\!\!\int (q(\varrho_\tau(y,z+ws)) - q(\varrho_\tau(y,z))) \\ &\qquad\qquad \times f(y)f'(z)\,dy\,dz\Big|. \end{aligned}$$

Using the Lipschitz property (4.12) of $\mathcal{Q}$, we see that $T_2(s)=O(s)$ and $T_1(s)=O(s)+O(T_3(s))$, where

$$T_3(s) = \sup_{0\le w\le 1}\sup_{|\tau-\vartheta|\le\Delta}\sup_{q\in\mathcal{Q}}\int |r_\tau(\mathbf{x}_{-1},y+ws) - r_\tau(\mathbf{x}_{-1},y)|(1+|y|)^\zeta|f'(y)|\,dy.$$

This shows that (2.11) holds if $T_3(s)=O(s)$.

To obtain that the class $\bar{\mathcal{H}}^*_\eta$ is Donsker, we will impose the following conditions (B1) and (B2) on $\mathcal{Q}$ and the class $\mathcal{R}$ defined by

$$\mathcal{R} = \{r_\tau(\mathbf{x}_{-1},\cdot + r_\vartheta(\mathbf{x})+s) : |\tau-\vartheta|\le\Delta, |s|\le \eta+\tilde{\Delta}\},$$

with $\tilde{\Delta} = \sup_{|\tau-\vartheta|\le\Delta}|r_\tau(\mathbf{x}) - r_\vartheta(\mathbf{x})|$. Note that the growth condition (4.2) implies that $\mathcal{R}$ has an envelope of the form $A_{\mathcal{R}}(1+|y|)$ for some constant $A_{\mathcal{R}}$.

(B1) For some integer $k$ and every $\delta>0$, there are $N=N_\delta=O(\delta^{-k})$ elements $q_1,\dots,q_N$ in $\mathcal{Q}$ such that $\mathcal{Q}$ is covered by the brackets $[q_i-\delta V_{\mathcal{Q}}, q_i+\delta V_{\mathcal{Q}}]$, $i=1,\dots,N$.

(B2) The class $\mathcal{R}$ has $L_2(\mu)$-bracketing numbers of polynomial growth for $\mu(dy)=(1+|y|)^{2\zeta}f(y)\,dy$: For some integer $j$,

$$N_{[\cdot]}(\delta,\mathcal{R},L_2(\mu)) = O(\delta^{-j}).$$

These properties, the growth condition (4.2) and the Lipschitz property (4.12) of $\mathcal{Q}$ imply that the class $\mathcal{G}=\{(y,z)\mapsto q(y+s,z+t+a(y)) : a\in\mathcal{R}; |s|,|t|\le C\}$ has $L_2(F\times F)$-bracketing numbers with polynomial growth for each finite $C$. Indeed, for $C\ge A_{\mathcal{R}}$, we can consider brackets of the form

$$\begin{aligned} &q(y+u,z+v+\bar{a}(y)) \pm C_*(2\delta+|a^*(y)-a_*(y)|)w(y)w(z) \\ &\qquad \pm \delta V_{\mathcal{Q}}(y+u,z+v+\bar{a}(y)), \end{aligned}$$

where $w(x)=(1+|x|)^\zeta$, $u$ and $v$ belong to the grid $\{i\delta : i\in\mathbb{Z}, |i\delta|\le B\}$ and $\bar{a}$ is the midpoint of a bracket $[a_*,a^*]$ for $\mathcal{R}$. Since $\mathcal{G}$ has polynomial growth, so do the classes $\mathcal{G}_1=\{\int g(\cdot,z)f(z)\,dz : g\in\mathcal{G}\}$ and $\mathcal{G}_2=\{\int g(y,\cdot)f(y)\,dy : g\in$ $\mathcal{G}\}$. Hence, these classes are Donsker. Since subsets and sums of Donsker classes are Donsker, and since $\bar{\mathcal{H}}_\eta^* \subset \mathcal{G}_1 + \mathcal{G}_2$ for large enough $C$, we see that $\bar{\mathcal{H}}_\eta^*$ is Donsker. Thus, we have the following result.

THEOREM 4.1. *Suppose Assumptions* B, K *and* R *hold. Suppose the class* $\mathcal{Q}$ *has envelope* $V_{\mathcal{Q}}$ *given by* (4.1) *and satisfies the Lipschitz property* (4.12) *and the growth property* (B1). *Let* $f$ *have finite Fisher information for location and a finite absolute moment of order greater than* $2\gamma_1 + 2\gamma_2 + 1$. *Let* $\mathcal{R}$ *satisfy* (B2) *and let the autoregression functions satisfy the growth conditions* (4.2), *the differentiability condition* (4.4) *and*

$$\sup_{|\tau-\vartheta|\le\Delta} \sup_{q\in\mathcal{Q}} \int |r_\tau(\mathbf{x}_{-1}, y+s) - r_\tau(\mathbf{x}_{-1}, y)|(1+|y|)^\zeta |f'(y)|\, dy = O(s).$$

*Then*

$$\sup_{q\in\mathcal{Q}} \left| \tilde{\nu}(q) - \nu(\vartheta, q) - \frac{1}{n}\sum_{i=1}^n \bar{h}_{\vartheta,q}^*(\varepsilon_i) - [D_{\vartheta,q}^*]^\top (\tilde{\vartheta} - \vartheta) \right| = o_p(n^{-1/2});$$

$$\sup_{q\in\mathcal{Q}} \left| \tilde{\nu}_w(q) - \nu(\vartheta, q) - \frac{1}{n}\sum_{i=1}^n \bar{h}_{\vartheta,q}^\#(\varepsilon_i) - [D_{\vartheta,q}^\#]^\top (\tilde{\vartheta} - \vartheta) \right| = o_p(n^{-1/2}),$$

*where* $D_{\vartheta,q}^* = \Psi_{\vartheta,q} - B(\bar{h}_{\vartheta,q}^*)$ *and* $D_{\vartheta,q}^\# = \Psi_{\vartheta,q} - B(\bar{h}_{\vartheta,q}^\#)$, *with* $\Psi_{\vartheta,q}$ *as given in Lemma* 4.1.

Exactly as in Sections 2 and 3, one obtains functional central limit theorems for the von Mises statistics and empirical estimators for their asymptotic covariance functions.

**5. Conditional distribution functions and quantiles of lag two.** The conditional distribution function $F_{\mathbf{x}}$ of the pair $(X_{n+1}, X_{n+2})$, given $\mathbf{X}_n = \mathbf{x}$, is defined by

$$\begin{aligned} F_{\mathbf{x}}(t,u) &= P(\varepsilon_1 + r_\vartheta(\mathbf{x}) \le t, \varepsilon_2 + r_\vartheta(\mathbf{x}_{-1}, \varepsilon_1 + r_\vartheta(\mathbf{x})) \le u) \\ &= \int_{-\infty}^t F(u - r_\vartheta(\mathbf{x}_{-1}, y)) f(y - r_\vartheta(\mathbf{x}))\, dy. \end{aligned}$$

This can also be written as

$$F_{\mathbf{x}}(t,u) = \iint q_{t,u}(\varrho_\vartheta(y,z)) f(y) f(z)\, dy\, dz$$

with $q_{t,u}(v,w) = \mathbf{1}[v \le t, w \le u]$. We estimate $F_{\mathbf{x}}(t,u)$ by

$$\tilde{F}_{\mathbf{x}}(t,u) = \iint q_{t,u}(\varrho_{\tilde{\vartheta}}(y,z)) \tilde{f}(y) \tilde{f}(z)\, dy\, dz$$

and its weighted version

$$\tilde{F}_{\mathbf{x}w}(t,u) = \iint q_{t,u}(\varrho_{\tilde{\vartheta}}(y,z))\tilde{f}_w(y)\tilde{f}_w(z)\,dy\,dz.$$

Here the class $\mathcal{Q}$ equals $\{q_{t,u} : t,u \in \mathbb{R}\}$. It has envelope $V_{\mathcal{Q}} = 1$; thus, condition (4.1) holds with $\gamma_1 = \gamma_2 = 0$ and $C_{\mathcal{Q}} = 1$. We have

$$\bar{h}^{(1)}_{\tau,q_{t,u}}(y) = F(u - r_\tau(\mathbf{x}_{-1}, y + r_\tau(\mathbf{x})))\mathbf{1}[y \le t - r_\tau(\mathbf{x})],$$

$$\bar{h}^{(2)}_{\tau,q_{t,u}}(z) = \int_{-\infty}^{t} \mathbf{1}[z \le u - r_\tau(\mathbf{x}_{-1}, y)] f(y - r_\tau(\mathbf{x}))\,dy.$$

We shall now show that $\bar{\mathcal{H}}^*_\eta$ is Donsker if the class $\mathcal{R}$ has $L_2(F)$-bracketing numbers with polynomial growth:

$$N_{[\cdot]}(\delta, \mathcal{R}, L_2(F)) = O(\delta^{-j}) \qquad \text{for some positive integer } j. \tag{5.1}$$

It is easy to check that $\bar{\mathcal{H}}^*_\eta \subset \mathcal{F}_1 + \mathcal{F}_2$, where

$$\mathcal{F}_1 = \{F(u - a(\cdot))\mathbf{1}[\cdot \le v] : a \in \mathcal{R}; u, v \in \mathbb{R}\},$$

$$\mathcal{F}_2 = \left\{\int_{-\infty}^{v} \mathbf{1}[\cdot \le u - a(y)] f(y)\,dy : a \in \mathcal{R}; u, v \in \mathbb{R}\right\}.$$

Since subsets and sums of Donsker classes are Donsker, it suffices to show that $\mathcal{F}_1$ and $\mathcal{F}_2$ are Donsker classes. For this, it is enough to show the classes $\mathcal{F}_1$ and $\mathcal{F}_2$ have $L_2(F)$-bracketing numbers with polynomial growth. For $\mathcal{F}_1$, take brackets of the form $[b_*, b^*] = [F(u_* - a^*(\cdot))\mathbf{1}[\cdot \le v_*], F(u^* - a_*(\cdot))\mathbf{1}[\cdot \le v^*]]$, where $[a_*, a^*]$ is an $(\varepsilon/\|f\|_\infty)$-bracket for $\mathcal{R}$; $v_*$, $v^*$ are chosen such that $F(v^*) - F(v_*) \le \varepsilon^2$; and $u_*$, $u^*$ are chosen such that either $u^* - u_* \le \varepsilon/\|f\|_\infty$, or $u_* = -\infty$ and (i) $\int F^2(u^* + A_{\mathcal{R}}(1+|y|))f(y)\,dy \le \varepsilon^2$, or $u^* = \infty$ and (ii) $\int (1 - F(u_* - A_{\mathcal{R}}(1+|y|)))^2 f(y)\,dy \le \varepsilon^2$. Then $[b_*, b^*]$ is a $3\varepsilon$-bracket for $\mathcal{F}_1$. Since $F$ has finite second moment, $t^2F(t)(1-F(t)) \to 0$ as $|t| \to \infty$. Using this, it is easy to see that $u^*$ in (i) can be chosen proportional to $-1/\varepsilon$ and $u_*$ in (ii) can be taken proportional to $1/\varepsilon$. Thus, under (5.1), we can cover $\mathcal{F}_1$ with $O(\varepsilon^{-j-4})$ brackets of this form.

For $\mathcal{F}_2$, take brackets of the form

$$[b_*, b^*] = \left[\int_{-\infty}^{v_*} \mathbf{1}[z \le u_* - a^*(y)] f(y)\,dy, \int_{-\infty}^{v^*} \mathbf{1}[z \le u^* - a_*(y)] f(y)\,dy\right],$$

where $[a_*, a^*]$ is an $\varepsilon^2/\|f\|_\infty$-bracket for $\mathcal{R}$; $F(v^*) - F(v_*) \le \varepsilon^2$; and $v_* \le v^*$ are chosen such that either $v^* - v_* \le \varepsilon^2/\|f\|_\infty$, or $v_* = -\infty$ and (i) $\int F(v^* + A_{\mathcal{R}}(1+|y|))f(y)\,dy \le \varepsilon^2$, or $v^* = \infty$ and (ii) $\int (1 - F(v_* - A_{\mathcal{R}}(1+|y|)))^2 f(y)\,dy \le \varepsilon^2$. Then $[b_*, b^*]$ is a $3\varepsilon$-bracket for $\mathcal{F}_2$. It is easy to check that, under (5.1), we can cover $\mathcal{F}_2$ with $O(\varepsilon^{-2j-5})$ brackets of this form.

This shows that condition (5.1) implies that $\bar{\mathcal{H}}^*_\eta$ is Donsker. In view of Remark 4.1, we then obtain that condition (4.8) is met. Using the moment inequality and interchanging the order of integration, we can bound the left-hand side of condition (4.9) by

$$\sup_{u\in\mathbb{R}}\int |F(u-r_\tau(\mathbf{x}_{-1},y))-F(u-r_\vartheta(\mathbf{x}_{-1},y))|f(y-r_\vartheta(\mathbf{x}))\,dy$$
$$\le \|f\|_\infty \int |r_\tau(\mathbf{x}_{-1},y)-r_\vartheta(\mathbf{x}_{-1},y)|f(y-r_\vartheta(\mathbf{x}))\,dy.$$

Thus, we have (4.9) in view of the Lebesgue dominated convergence theorem and the growth condition (4.2). Lemmas 4.1 and 4.2 now imply conditions (2.9) and (2.10) of Theorem 2.2. Finally, (2.11) is implied by the two conditions

$$\sup_{\tau,t,u}\left|\int_{-\infty}^{t} F(u-r_\tau(\mathbf{x}_{-1},y+r_\tau(\mathbf{x})))(f(y-s)-f(y)+sf'(y))\,dy\right| = O(s^2) \tag{5.2}$$

and

$$\sup_{\tau,t,u}\left|\int_{-\infty}^{t}(F(u-r_\tau(\mathbf{x}_{-1},y)-s) - F(u-r_\tau(\mathbf{x}_{-1},y))+sf(u-r_\tau(\mathbf{x}_{-1},y)))f(y-r_\tau(\mathbf{x}))\,dy\right| = O(s^2),$$

where the suprema extend over all real $t$ and $u$ and all $\tau$ with $|\tau-\vartheta|\le\Delta$. The latter condition is satisfied if $f$ is Lipschitz. If we set $a_\tau(y)=r_\tau(\mathbf{x}_{-1},y+r_\tau(\mathbf{x}))$ and use (4.13), we obtain that the integral in (5.2) can be written

$$-s\int_0^1\left(\int_{-\infty}^{t-ws}F(u-a_\tau(y+ws))f'(y)\,dy-\int_{-\infty}^{t}F(u-a_\tau(y))f'(y)\,dy\right)dw.$$

If $f$ is Lipschitz, so that $f'$ is bounded, we see that condition (5.2) is implied by

$$\sup_{|\tau-\vartheta|\le\Delta}\int |r_\tau(\mathbf{x}_{-1},y+s)-r_\tau(\mathbf{x}_{-1},y)||f'(y-r_\tau(\mathbf{x}))|\,dy = O(s). \tag{5.3}$$

If $f$ has finite Fisher information, then $f'$ is integrable and a sufficient condition for (5.3) is that there is a constant $L$ such that

$$|r_\tau(\mathbf{x}_{-1},y_1)-r_\tau(\mathbf{x}_{-1},y_2)|\le L|y_1-y_2|, \qquad y_1,y_2\in\mathbb{R};|\tau-\vartheta|\le\Delta. \tag{5.4}$$

Hence, Theorem 2.2 gives the following stochastic expansions for $\tilde F_{\mathbf{x}}$ and $\tilde F_{\mathbf{x}w}$.

THEOREM 5.1. *Let Assumptions B, F, K and R hold and let $f$ be Lipschitz. Suppose that (4.4), the growth conditions (4.2) and (5.1), and (5.3) or (5.4) hold. Then*

$$\sup_{t,u\in\mathbb{R}}\left|\tilde{F}_{\mathbf{x}}(t,u)-F_{\mathbf{x}}(t,u)-\frac{1}{n}\sum_{i=1}^{n}h^*_{t,u}(\varepsilon_i)-[D^*_{t,u}]^\top(\tilde{\vartheta}-\vartheta)\right|=o_p(n^{-1/2});$$

$$\sup_{t,u\in\mathbb{R}}\left|\tilde{F}_{\mathbf{x}w}(t,u)-F_{\mathbf{x}}(t,u)-\frac{1}{n}\sum_{i=1}^{n}h^{\#}_{t,u}(\varepsilon_i)-[D^{\#}_{t,u}]^\top(\tilde{\vartheta}-\vartheta)\right|=o_p(n^{-1/2}),$$

*where $D^*_{t,u}=\Psi_{t,u}-B(h^*_{t,u})$ and $D^{\#}_{t,u}=\Psi_{t,u}-B(h^{\#}_{t,u})$ and where*

$$\begin{aligned}h_{t,u}(y)&=F(u-r_\vartheta(\mathbf{x}_{-1},y+r_\vartheta(\mathbf{x})))\mathbf{1}[y\le t-r_\vartheta(\mathbf{x})]\\&\qquad+\int_{-\infty}^{t}\mathbf{1}[y\le u-r_\vartheta(\mathbf{x}_{-1},z)]f(z-r_\vartheta(\mathbf{x}))\,dz;\\\Psi_{t,u}&=-\int_{-\infty}^{t}F(u-r_\vartheta(\mathbf{x}_{-1},y))f'(y-r_\vartheta(\mathbf{x}))\,dy\,\dot{r}_\vartheta(\mathbf{x})\\&\qquad-\int_{-\infty}^{t}f(u-r_\vartheta(\mathbf{x}_{-1},y))f(y-r_\vartheta(\mathbf{x}))\dot{r}_\vartheta(\mathbf{x}_{-1},y)\,dy.\end{aligned}$$

EXAMPLE 5.1. Consider the AR(1) model, in which $r_\vartheta(x)=\vartheta x$ and $|\vartheta|<1$. Clearly, Assumption R and conditions (4.2), (4.4) and (5.1) hold. Also, condition (5.4) holds with $L=1$. Thus, if $f$ has finite Fisher information for location and is Lipschitz, then all the assumptions of Theorem 5.1 can be met.

EXAMPLE 5.2. Consider the EXPAR(1) model, in which $\vartheta=(\vartheta_1,\vartheta_2)$ with $\vartheta_1<1$, and $r_\vartheta(x)=(\vartheta_1+\vartheta_2\exp(-\gamma x^2))x$. Here the exponent $\gamma$ is assumed known. The assumptions of Theorem 5.1 can be met if $f$ has finite Fisher information for location and is Lipschitz. Clearly, Assumption R and conditions (4.2), (4.4) and (5.1) hold. Moreover, condition (5.4) is satisfied.

We have phrased the conditions on $r_\tau$ in Theorem 5.1 sufficiently general to cover discontinuous autoregression functions such as those appearing in SETAR models.

EXAMPLE 5.3. Consider the SETAR(2, 1, 1) model with known threshold $\xi$. In this model, $\vartheta=(\vartheta_1,\vartheta_2)$ with $\vartheta_1<1$, $\vartheta_2<1$, $\vartheta_1\vartheta_2<1$, and $r_\vartheta(x)=\vartheta_1x\mathbf{1}[x\le\xi]+\vartheta_2x\mathbf{1}[x>\xi]$. It is easily seen that Assumption R and the conditions (4.2), (4.4) and (5.1) hold. Suppose now that $f$ has finite Fisher information for location and is Lipschitz. If $\xi=0$, then the Lipschitz condition (5.4) holds. If $\xi\ne0$, then the Lipschitz condition (5.4) does not hold,

but one has

$$|r_\tau(y+s) - r_\tau(y)| \le (|\tau_1| + |\tau_2|)(|s| + \mathbf{1}[|y-\xi| \le |s|]).$$

This and the fact that $f'$ is bounded and integrable yield (5.3).

The one-dimensional lag-two conditional distribution function $G_{\mathbf{x}}(u) = F_{\mathbf{x}}(\infty, u)$ at $u$ of $X_{n+2}$, given $\mathbf{X}_n = \mathbf{x}$, is

$$G_{\mathbf{x}}(u) = \int F(u - r_\vartheta(\mathbf{x}_{-1}, y)) f(y - r_\vartheta(\mathbf{x}))\, dy.$$

We estimate $G_{\mathbf{x}}(u)$ by

$$\tilde{G}_{\mathbf{x}}(u) = \iint \mathbf{1}[z + r_{\tilde\vartheta}(\mathbf{x}_{-1}, y + r_{\tilde\vartheta}(\mathbf{x})) \le u] \tilde{f}(y)\tilde{f}(z)\, dy\, dz$$

and its weighted version

$$\tilde{G}_{\mathbf{x}w}(u) = \iint \mathbf{1}[z + r_{\tilde\vartheta}(\mathbf{x}_{-1}, y + r_{\tilde\vartheta}(\mathbf{x})) \le u] \tilde{f}_w(y)\tilde{f}_w(z)\, dy\, dz.$$

We obtain stochastic expansions as in Theorem 5.1, with $t$ replaced by $\infty$.

THEOREM 5.2. *Let Assumptions* B, F, K *and* R *hold and let $f$ be Lipschitz. Suppose that* (4.4), *the growth conditions* (4.2) *and* (5.1) *and* (5.3) *hold. Then*

$$\sup_{u\in\mathbb{R}} \left| \tilde{G}_{\mathbf{x}}(u) - G_{\mathbf{x}}(u) - \frac{1}{n}\sum_{i=1}^n h_u^*(\varepsilon_i) - [D_u^*]^\top (\tilde\vartheta - \vartheta) \right| = o_p(n^{-1/2});$$

$$\sup_{u\in\mathbb{R}} \left| \tilde{G}_{\mathbf{x}w}(u) - G_{\mathbf{x}}(u) - \frac{1}{n}\sum_{i=1}^n h_u^\#(\varepsilon_i) - [D_u^\#]^\top (\tilde\vartheta - \vartheta) \right| = o_p(n^{-1/2}),$$

*where* $D_u^* = \Psi_u - B(h_u^*)$ *and* $D_u^\# = \Psi_u - B(h_u^\#)$ *and where*

$$h_u(y) = F(u - r_\vartheta(\mathbf{x}_{-1}, y + r_\vartheta(\mathbf{x}))) + \int \mathbf{1}[y \le u - r_\vartheta(\mathbf{x}_{-1}, z)] f(z - r_\vartheta(\mathbf{x}))\, dz;$$

$$\Psi_u = -\int F(u - r_\vartheta(\mathbf{x}_{-1}, y)) f'(y - r_\vartheta(\mathbf{x}))\, dy\ \dot{r}_\vartheta(\mathbf{x})$$
$$\qquad - \int f(u - r_\vartheta(\mathbf{x}_{-1}, y)) f(y - r_\vartheta(\mathbf{x})) \dot{r}_\vartheta(\mathbf{x}_{-1}, y)\, dy.$$

EXAMPLE 5.4. We apply Theorem 5.2 to the conditional quantile function of lag two. The conditional $v$-quantile of $X_{n+2}$, given $\mathbf{X}_n = \mathbf{x}$, is the left-inverse $G_{\mathbf{x}}^{-1}(v)$ of $G_{\mathbf{x}}$ at $v$. We estimate it by $\tilde{G}_{\mathbf{x}w}^{-1}(v)$. Since $f$ was assumed positive, $G_{\mathbf{x}}$ has a positive density

$$g_{\mathbf{x}}(u) = \int f(u - r_\vartheta(\mathbf{x}_{-1}, y)) f(y - r_\vartheta(\mathbf{x}))\, dy.$$

Let $0 < c \le d < 1$. As in Example 3.1, we use Proposition 1 of [10] to obtain the stochastic expansion

$$\sup_{v\in[c,d]}\left|\tilde{G}_{\mathbf{x}w}^{-1}(v)-G_{\mathbf{x}}^{-1}(v)+\frac{1}{n}\sum_{i=1}^{n}\frac{h^{\#}_{G_{\mathbf{x}}^{-1}(v)}(\varepsilon_i)+[D^{\#}_{G_{\mathbf{x}}^{-1}(v)}]^{\top}(\tilde{\vartheta}-\vartheta)}{g_{\mathbf{x}}(G_{\mathbf{x}}^{-1}(v))}\right|=o_p(n^{-1/2}).$$

It follows that the smoothed and weighted lag-two conditional quantile process

$$\{n^{1/2}(\tilde{G}_{\mathbf{x}w}^{-1}(v)-G_{\mathbf{x}}^{-1}(v)) : v\in[c,d]\}$$

converges weakly in $\ell_\infty([c,d])$ to a centered Gaussian process.

EXAMPLE 5.5. Consider the AR(1) model $X_i=\vartheta X_{i-1}+\varepsilon_i$ with $|\vartheta|<1$. Let us also take $\tilde{\vartheta}$ to be the sample correlation coefficient, which satisfies (3.2). We are interested in predicting the probability $G_x(u)=P(X_{n+2}\le u|X_n=x)$, which can be expressed as

$$G_x(u)=P(\varepsilon_2+\vartheta\varepsilon_1+\vartheta^2x\le u)=\int F(u-\vartheta y-\vartheta^2x)f(y)\,dy.$$

We assume that $f$ has finite Fisher information for location and is Lipschitz, so that the requirements of Theorem 5.1 and, hence, of Theorem 5.2 are met as demonstrated in Example 5.1. The smoothed von Mises estimator is

$$\tilde{G}_x(u)=\tilde{F}_{\mathbf{x}}(\infty,u)=\iint\mathbf{1}[z+\tilde{\vartheta}y+\tilde{\vartheta}^2x\le u]\hat{f}(y)\hat{f}(z)\,dy\,dz,$$

and its weighted counterpart is

$$\tilde{G}_{xw}(u)=\tilde{F}_{\mathbf{x}w}(\infty,u)=\iint\mathbf{1}[z+\tilde{\vartheta}y+\tilde{\vartheta}^2x\le u]\hat{f}_w(y)\hat{f}_w(z)\,dy\,dz.$$

Since $\dot{r}_\vartheta(x)=x$ and $E[X]=0$, we see that $B(g)=0$ for all $g\in L_2(F)$. Thus, we obtain from Theorem 5.2 and from expansion (3.2) for $\tilde{\vartheta}$ that

$$\tilde{G}_x(u)=G_x(u)+\frac{1}{n}\sum_{i=1}^{n}\left(h_u(\varepsilon_i)-2G_x(u)+\Psi_u\frac{1-\vartheta^2}{\sigma^2}X_{i-1}\varepsilon_i\right)+o_p(n^{-1/2})$$

and

$$\tilde{G}_{xw}(u)=\tilde{G}_x(u)-\frac{c_u}{\sigma^2}\frac{1}{n}\sum_{i=1}^{n}\varepsilon_i+o_p(n^{-1/2}),$$

where $\Psi_u=-E[(\varepsilon+2\vartheta x)f(u-\vartheta\varepsilon-\vartheta^2x)]$, $c_u=E[\varepsilon h_u(\varepsilon)]$ and

$$h_u(\varepsilon)=h_{\infty,u}(\varepsilon)=\begin{cases}F(u-\vartheta\varepsilon-\vartheta^2x)+F((u-\varepsilon-\vartheta^2x)/\vartheta), & \vartheta>0,\\ F(u)+\mathbf{1}[\varepsilon\le u], & \vartheta=0,\\ F(u-\vartheta\varepsilon-\vartheta^2x)+1-F((u-\varepsilon-\vartheta^2x)/\vartheta), & \vartheta<0.\end{cases}$$

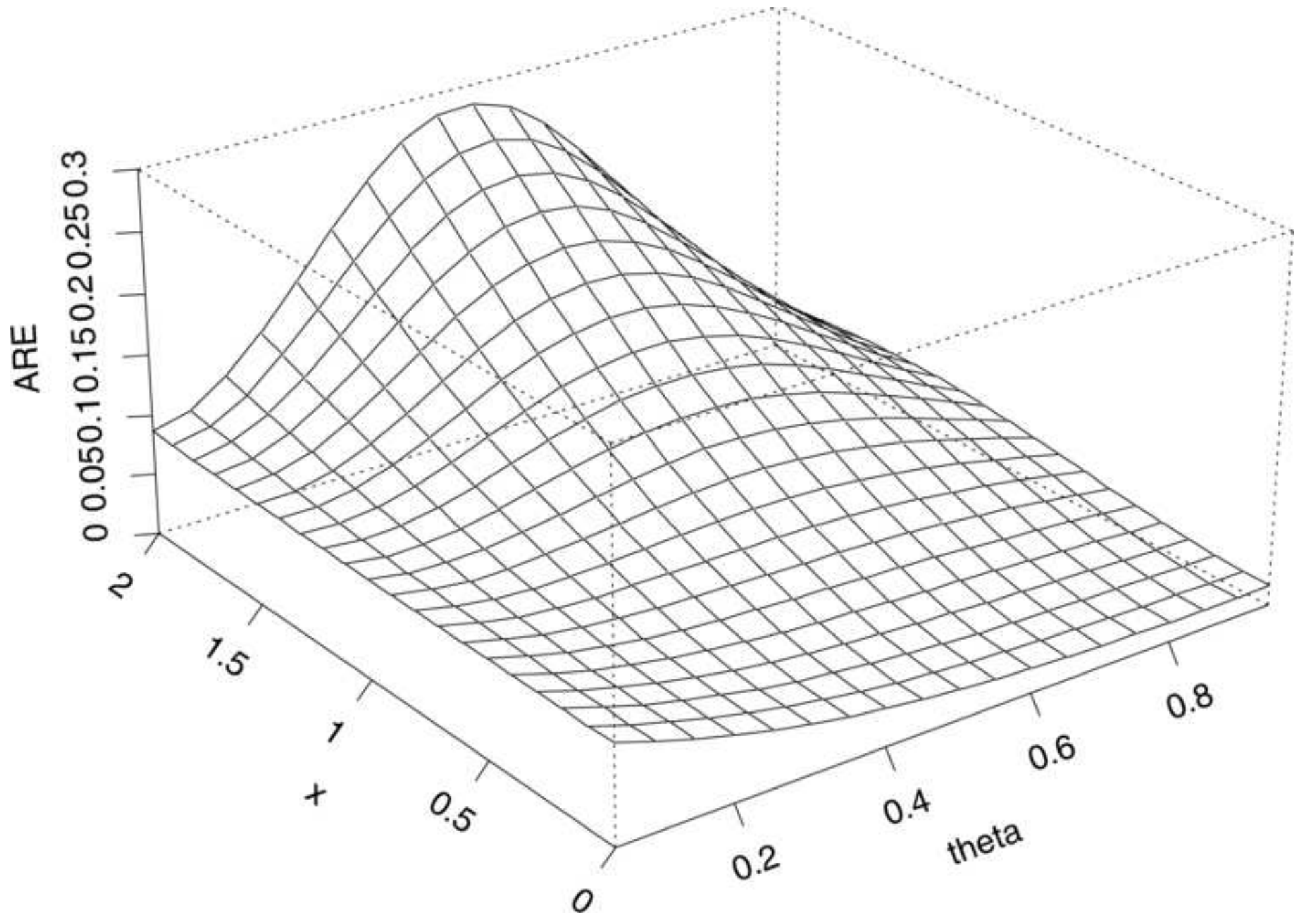

FIG. 1. *The asymptotic relative efficiency* $\tau_w^2/\tau^2$ *of the unweighted versus the weighted estimator for* $u = 0$ *for various values of* $\vartheta$ *and* $x$.

Consequently, $n^{1/2}(\tilde{G}_x(u) - G_x(u))$ is asymptotically normal with mean zero and variance $\tau^2 = \mathrm{Var}(h_u(\varepsilon)) + \Psi_u^2(1 - \vartheta^2)$, while $n^{1/2}(\tilde{G}_{xw}(u) - G_x(u))$ is asymptotically normal with mean zero and variance $\tau_w^2 = \tau^2 - c_u^2/\sigma^2$. Therefore, the weighted version has a smaller asymptotic variance unless $c_u = 0$. The variance reductions can be considerable. Figure 1 is a graph of the asymptotic relative efficiency $\tau_w^2/\tau^2$ of the unweighted with respect to the weighted estimator as a function of $\vartheta$ (ranging from 0.05 to 0.95) and $x$ (ranging from 0 to 2) in the case of the standard normal density $f$ and $u = 0$. As one can see from the graph, the ratio is always below 0.3 and can be as small as 0.0151. Thus, variance reductions of over 98% are possible.

**6. Efficiency.** In this section we prove that the weighted versions of our estimators are efficient. We recall that, among all "regular" estimators, an estimator for a vector-valued functional is *efficient* in the sense of Hájek and Le Cam if its standardized error is asymptotically maximally concentrated in symmetric convex sets. In a locally asymptotically normal model, an estimator for a differentiable functional is regular and efficient if and only if it is asymptotically linear with influence function equal to the canonical gradient of the functional. For our nonlinear autoregressive model, these concepts and the explicit form of the characterization are given in [32] for differentiable functionals of both the autoregression parameter and the innovation density.

Here we are interested in estimating the functional

$$\kappa(\vartheta, f) = \psi(h_\vartheta, f) = \int h_\vartheta(\mathbf{y}) \prod_{j=1}^m f(y_j)\, d\mathbf{y},$$

where $\{h_\tau : |\tau - \vartheta| \le \Delta\}$ is a class of measurable functions from $\mathbb{R}^m$ into $\mathbb{R}$. We assume that the class has an envelope $H$ of the form $H(y_1, \dots, y_m) = V(y_1) \cdots V(y_m)$ with $V$ satisfying Assumption V. We use $\bar{h}_\tau$ as defined in Section 2 and assume that

$$\int (\bar{h}_\tau - \bar{h}_\vartheta)^2\, dF \to 0 \qquad \text{as } \tau \to \vartheta. \tag{6.1}$$

LEMMA 6.1. *Suppose, in addition to the above, that $\tau \mapsto \psi(h_\tau, f)$ is differentiable at $\vartheta$ with gradient $\Psi_\vartheta$. Let $\vartheta_n$ be a sequence in $\mathbb{R}^d$ such that $n^{1/2}(\vartheta_n - \vartheta) \to u$. Let $f_n$ be a sequence of densities such that $\|f_n - f\|_{V^2} \to 0$ and*

$$\int (n^{1/2}(f_n^{1/2}(y) - f^{1/2}(y)) - \tfrac{1}{2} v(y) f^{1/2}(y))^2\, dy \to 0$$

*for some $v \in L_2(F)$. Then*

$$n^{1/2}(\psi(h_{\vartheta_n}, f_n) - \psi(h_\vartheta, f)) \to \Psi_\vartheta^\top u + \int \bar{h}_\vartheta v\, dF.$$

PROOF. Express $\psi(h_{\vartheta_n}, f_n) - \psi(h_\vartheta, f)$ as the sum $T_1 + T_2 + T_3$ with

$$T_1 = \psi(h_{\vartheta_n}, f_n) - \psi(h_{\vartheta_n}, f) - \int \bar{h}_{\vartheta_n}(y)(f_n(y) - f(y))\, dy,$$

$$T_2 = \int \bar{h}_{\vartheta_n}(y)(f_n(y) - f(y))\, dy,$$

$$T_3 = \psi(h_{\vartheta_n}, f) - \psi(h_\vartheta, f).$$

We have $n^{1/2} T_3 \to \Psi_\vartheta^\top u$. The argument given in the proof of Theorem 2.1 shows that $T_1 = O(\|f_n - f\|_V^2)$. Writing $s_n = f_n^{1/2}$, $s = f^{1/2}$ and $f_n - f = (s_n - s) \times (s_n + s)$, and applying the Cauchy–Schwarz inequality, we obtain that $\|f_n - f\|_V^2 \le 2(\|f_n\|_{V^2} + \|f\|_{V^2}) \|(s_n - s)^2\|_1 = O(n^{-1})$. Thus, $n^{1/2} T_1 \to 0$. Finally, $n^{1/2} T_2 \to \int \bar{h}_\vartheta v\, dF$ by the same argument as for Lemma 7.2. □

As in Section 2, let $\mathcal{V}$ be the set of all $v \in L_2(F)$ with $\int v(y)\, dF(y) = 0$ and $\int y v(y)\, dF(y) = 0$. For each $v \in \mathcal{V}$, there is a sequence $f_n = f_{nv}$ of zero mean densities as required in the previous lemma. As shown in [32], these densities can be chosen to also satisfy $\|(f_n - f)/f\|_\infty \to 0$ and to have finite Fisher information for location if $f$ has it. We also require now that (3.1) in [32] hold: for every sequence $\vartheta_n$ and $f_n = f_{nv}$ as above, the corresponding stationary density converges in $L_1$ to the density of $G$. Under this assumption, one has local asymptotic normality. As seen in Section 2, under appropriate conditions, the estimator $\psi(h_{\tilde\vartheta}, \tilde f_w)$ has the stochastic expansion

$$\psi(h_{\tilde\vartheta}, \tilde f_w) = \psi(h_\vartheta, f) + \frac{1}{n}\sum_{i=1}^n \bar h_\vartheta^{\#}(\varepsilon_i) + \left(\Psi_\vartheta - \mu \int \bar h_\vartheta^{\#} \ell \, dF\right)^\top (\tilde\vartheta - \vartheta) + o_p(n^{-1/2}),$$

with $\mu = E[\dot r_\vartheta(\mathbf{X})]$ and $\bar h_\vartheta^{\#}(y) = \bar h_\vartheta(y) - \int \bar h_\vartheta \, dF - \sigma^{-2} y \int u \bar h_\vartheta(u)\, dF(u)$. Recall that $\bar h_\vartheta^{\#}$ is the projection of $\bar h_\vartheta$ onto $\mathcal{V}$. The projection of $\ell$ onto $\mathcal{V}$ is $\ell^{\#}(y) = \ell(y) - \sigma^{-2} y$. If $\tilde\vartheta$ is efficient, then, by characterization (3.12) of [32], it has the stochastic expansion

$$\tilde\vartheta = \vartheta + \Lambda^{-1} \frac{1}{n}\sum_{i=1}^n S(\mathbf{X}_{i-1}, \varepsilon_i) + o_p(n^{-1/2}),$$

where $S(\mathbf{X}, \varepsilon) = \dot r_\vartheta(\mathbf{X})\ell(\varepsilon) - \mu \ell^{\#}(\varepsilon)$ and $\Lambda = E[S(\mathbf{X}, \varepsilon) S(\mathbf{X}, \varepsilon)^\top] = JR - J^{\#}\mu\mu^\top$, with $J$ and $J^{\#}$ the second moments of $\ell(\varepsilon)$ and $\ell^{\#}(\varepsilon)$, and $R = E[\dot r_\vartheta(\mathbf{X}) \dot r_\vartheta(\mathbf{X})^\top]$. If an efficient estimator $\tilde\vartheta$ is used in $\psi(h_{\tilde\vartheta}, \tilde f_w)$, we obtain the stochastic expansion

$$\psi(h_{\tilde\vartheta}, \tilde f_w) = \psi(h_\vartheta, f) + \frac{1}{n}\sum_{i=1}^n S_{\#}(\mathbf{X}_{i-1}, \varepsilon_i) + o_p(n^{-1/2}),$$

with

$$S_{\#}(\mathbf{X}, \varepsilon) = \bar h_\vartheta^{\#}(\varepsilon) + M^\top (\dot r_\vartheta(\mathbf{X})\ell(\varepsilon) - \mu \ell^{\#}(\varepsilon)),$$
$$M = \Lambda^{-1}\left(\Psi_\vartheta - \mu \int \bar h_\vartheta^{\#} \ell \, dF\right).$$

For $v \in \mathcal{V}$, we have

$$\begin{aligned} E[S_{\#}(\mathbf{X}, \varepsilon) v(\varepsilon)] &= \int \bar h_\vartheta^{\#} v \, dF + M^\top \mu \left(\int \ell v \, dF - \int \ell^{\#} v \, dF\right) \\ &= \int \bar h_\vartheta^{\#} v \, dF \\ &= \int \bar h_\vartheta v \, dF. \end{aligned}$$

Furthermore,

$$\begin{aligned}E[S_\#(\mathbf{X},\varepsilon)\dot r_\vartheta(\mathbf{X})\ell(\varepsilon)] &= \mu\int \bar h_\vartheta^\# \ell\, dF + (JR - J^\#\mu\mu^\top)M\\ &= \mu\int \bar h_\vartheta^\# \ell\, dF + \Lambda\Lambda^{-1}\Big(\Psi_\vartheta - \mu\int \bar h_\vartheta^\# \ell\, dF\Big)\\ &= \Psi_\vartheta.\end{aligned}$$

This shows that, for all $u \in \mathbb{R}^d$ and $v \in \mathcal{V}$,

$$E[S_\#(\mathbf{X},\varepsilon)(u^\top \dot r_\vartheta(\mathbf{X})\ell(\varepsilon) + v(\varepsilon))] = u^\top\Psi_\vartheta + \int \bar h_\vartheta v\, dF.$$

Since $S_\#(\mathbf{X},\varepsilon)$ is of the form $S_\#(\mathbf{X},\varepsilon) = u_0^\top \dot r_\vartheta(\mathbf{X})\ell(\varepsilon) + v_0(\varepsilon)$ for some $u_0 \in \mathbb{R}^d$ and $v_0 \in \mathcal{V}$, we obtain that $S_\#(\mathbf{X},\varepsilon)$ is the canonical gradient of the functional $\psi(h_\vartheta, f)$. Hence, $\psi(h_{\tilde\vartheta}, \tilde f_w)$ is efficient by the characterization (3.5) in [32], provided $S_\#(\mathbf{X},\varepsilon)$ is almost surely not zero.

The stochastic expansion of $\psi(h_{\tilde\vartheta}, \tilde f_w)$ given above implies that $n^{1/2} \times (\psi(h_{\tilde\vartheta}, \tilde f_w) - \psi(h_\vartheta, f))$ is asymptotically normal with variance

$$\begin{aligned}&E[\bar h_\vartheta^2(\varepsilon)] - (E[\bar h_\vartheta(\varepsilon)])^2 - \sigma^{-2}(E[\varepsilon\bar h_\vartheta(\varepsilon)])^2\\ &\quad + (\Psi_\vartheta - \mu E[\bar h_\vartheta(\varepsilon)\ell^\#(\varepsilon)])^\top \Lambda(\Psi_\vartheta - \mu E[\bar h_\vartheta(\varepsilon)\ell^\#(\varepsilon)]).\end{aligned}$$

## 7. Technical details.

We begin with a characterization of compact subsets of $L_2(\nu)$ for a measure $\nu$ with Lebesgue density.

LEMMA 7.1. *Let $\nu$ be a finite measure with Lebesgue density $\varphi$. Let $W \in L_2(\nu)$ satisfy*

$$\int (W(x-s) - W(x))^2\,\nu(dx) \to 0 \qquad \text{as } s \to 0. \tag{7.1}$$

*Then a subset $\mathcal{G}$ of $L_2(\nu)$ with envelope $W$ is totally bounded if and only if*

$$\sup_{g\in\mathcal{G}} \int (g(x-s) - g(x))^2\,\nu(dx) \to 0 \qquad \text{as } s \to 0. \tag{7.2}$$

PROOF. Let $\lambda$ denote the Lebesgue measure. Let $\bar{\mathcal{G}}$ denote the closure of $\mathcal{G}$ in $L_2(\nu)$. Clearly, $\mathcal{G}$ is totally bounded if and only if $\bar{\mathcal{G}}$ is compact in $L_2(\nu)$. The latter is equivalent to compactness of $\bar{\mathcal{G}}\sqrt{\varphi}$ in $L_2(\lambda)$. By the Fréchet–Kolmogorov theorem (see [45], page 275), compactness of $\bar{\mathcal{G}}\sqrt{\varphi}$ is equivalent to

$$\sup_{g\in\mathcal{G}} \int g^2\varphi\, d\lambda < \infty, \tag{7.3}$$

$$\sup_{g\in\mathcal{G}}\int (g(x-s)\sqrt{\varphi(x-s)}-g(x)\sqrt{\varphi(x)}\,)^2\,dx \to 0 \qquad \text{as } s\to 0, \tag{7.4}$$

$$\sup_{g\in\mathcal{G}}\int_{|x|>M} g^2(x)\varphi(x)\,dx \to 0 \qquad \text{as } M\to\infty. \tag{7.5}$$

Since $\mathcal{G}$ has envelope $W$ in $L_2(\nu)$, properties (7.3) and (7.5) are automatically satisfied. Since

$$\begin{aligned}((g(x-s)&\sqrt{\varphi(x-s)}-g(x)\sqrt{\varphi(x)}\,)-(g(x-s)-g(x))\sqrt{\varphi(x)}\,)^2\\ &= g^2(x-s)(\sqrt{\varphi(x-s)}-\sqrt{\varphi(x)}\,)^2,\end{aligned}$$

properties (7.4) and (7.2) are equivalent if

$$\int W^2(x-s)(\sqrt{\varphi(x-s)}-\sqrt{\varphi(x)}\,)^2\,dx\to 0 \qquad \text{as } s\to 0. \tag{7.6}$$

The above identity with $g=W$, together with continuity of translation in $L_2(\lambda)$, for which we refer to [30], Theorem 9.5, shows that (7.6) and (7.1) are equivalent. $\square$

The next lemma discusses uniform differentiability of integrals with respect to Hellinger differentiable densities.

LEMMA 7.2. *Let $\{p_\tau : |\tau-\vartheta|\le\Delta\}$ be a family of densities with respect to some measure $\nu$. Let $p_\tau$ be Hellinger differentiable at $\vartheta$ with Hellinger derivative $\chi$. Let $W$ be a nonnegative function such that*

$$\int W^2 p_\tau\,d\nu \to \int W^2 p_\vartheta\,d\nu \qquad \textit{as } \tau\to\vartheta. \tag{7.7}$$

*Then*

$$\sup_{|g|\le W}\left|\int g(p_\tau-p_\vartheta-\chi^\top(\tau-\vartheta)p_\vartheta)\,d\nu\right| = o(|\tau-\vartheta|). \tag{7.8}$$

*Moreover, if $\{g_\tau : |\tau-\vartheta|\le\Delta\}$ has envelope $W$ and $\int(g_\tau-g_\vartheta)^2 p_\vartheta\,d\nu\to 0$, then*

$$\begin{aligned}&\int g_\tau(p_\tau-p_\vartheta)\,d\nu\\ &\qquad = \int g_\vartheta\chi^\top p_\vartheta\,d\nu\,(\tau-\vartheta)+o(|\tau-\vartheta|).\end{aligned} \tag{7.9}$$

PROOF. Hellinger differentiability implies that $p_\tau\to p_\vartheta$ in $\nu$-measure. This and (7.7) yield $\int W^2|p_\tau-p_\vartheta|\,d\nu\to 0$. Let $s_\tau=p_\tau^{1/2}$ and $r_\tau=s_\tau-s_\vartheta-\frac{1}{2}\chi^\top(\tau-\vartheta)s_\vartheta$. Hellinger differentiability means that $\int r_\tau^2\,d\nu=o(|\tau-\vartheta|^2)$.

Since $p_\tau - p_\vartheta - \chi^\top(\tau-\vartheta)p_\vartheta = r_\tau(s_\tau+s_\vartheta) + \frac{1}{2}\chi^\top(\tau-\vartheta)s_\vartheta(s_\tau-s_\vartheta)$, an application of the Cauchy–Schwarz inequality shows that the square of the left-hand side of (7.8) can be bounded by

$$2\int W^2(s_\tau+s_\vartheta)^2\,d\nu\int r_\tau^2\,d\nu$$
$$+\tfrac{1}{2}\int(\chi^\top(\tau-\vartheta))^2 p_\vartheta\,d\nu\int W^2(s_\tau-s_\vartheta)^2\,d\nu.$$

Using $(s_\tau+s_\vartheta)^2\le 2(p_\tau+p_\vartheta)$ and $(s_\tau-s_\vartheta)^2\le|p_\tau-p_\vartheta|$, we obtain (7.8). To prove (7.9), it therefore remains to show that $\int(g_\tau-g_\vartheta)\chi p_\vartheta\,d\nu\to 0$ as $\tau\to\vartheta$. But this is an easy consequence of $\int(g_\tau-g_\vartheta)^2 p_\vartheta\,d\nu\to 0$. $\square$

**Acknowledgments.** We thank two referees and an Associate Editor for many suggestions that have improved the presentation.

U. U. Müller
Department of Statistics
Texas A&M University
College Station, Texas 77843-3143
USA
E-mail: uschi@stat.tamu.edu

A. Schick
Department of Mathematical Sciences
Binghamton University
Binghamton, New York 13902-6000
USA
E-mail: anton@math.binghamton.edu

W. Wefelmeyer
Universität zu Köln
Mathematisches Institut
Weyertal 86-90
50931 Köln
Germany
E-mail: wefelm@math.uni-koeln.de